\documentclass[3p,numbers,sort&compress,nopreprintline]{elsarticle}
\nonstopmode
\usepackage{lineno,hyperref,amsmath,amssymb,bm, physics, filecontents, mathrsfs, xcolor}
\usepackage[mathscr]{euscript}
\date{}

\usepackage{moreverb}
\usepackage{float}
\usepackage{graphicx}
\renewcommand{\vec}[1]{\ensuremath{\boldsymbol{#1}}}

\newcommand\BibTeX{{\rmfamily B\kern-.05em \textsc{i\kern-.025em b}\kern-.08em
T\kern-.1667em\lower.7ex\hbox{E}\kern-.125emX}}

\usepackage{filecontents}

\usepackage{numcompress}

\begin{document}

\begin{frontmatter}

\title{A computational periporomechanics model for localized failure in unsaturated porous media}
\author{Shashank Menon}
\author{Xiaoyu Song \corref{mycorrespondingauthor}}
\address{Engineering School of Sustainable Infrastructure and Environment, University of Florida, Gainesville, FL-32611}

\cortext[mycorrespondingauthor]{Corresponding author}
\ead{xysong@ufl.edu}

\begin{abstract}
We implement a computational periporomechanics model for simulating localized failure in unsaturated porous media. The coupled periporomechanics model is based on the peridynamic state concept and the effective force state concept. The coupled governing equations are integral-differential equations without assuming the continuity of solid displacement and fluid pressures. The fluid flow and effective force states are determined by nonlocal fluid pressure and deformation gradients through the multiphase constitutive correspondence principle. The coupled periporomechanics model is numerically implemented for high-performance computing by an implicit multiphase meshfree method utilizing the message passing interface. The numerical implementation is validated by simulating classical poromechanics problems and comparing numerical results with analytical solutions and experimental data. Numerical examples are presented to demonstrate the robustness of the implemented computational periporomechanics model for simulating localized failure in unsaturated porous media.
\end{abstract}

\begin{keyword}
Periporomechanics \sep coupled \sep localized failure \sep unsaturated porous media \sep nonlocal \sep parallel computing
\end{keyword}

\end{frontmatter}

\section{Introduction}\label{sec1}

Strain localization of geomaterials is a significant problem because it is closely related to the onset of failure of geomaterials (e.g., \citep{rudnicki1975conditions,vardoulakis1995bifurcation, borja2002bifurcation, borja2013shear, song2017transient}). For example, the failure of soil slopes is a typical strain localization problem in which deformation usually occurs in a narrow banded zone with a finite thickness \citep{borja2004computationalI, lu2013hillslope}. Physically, strain localization in solid materials is characterized by the concentration of strain or displacement discontinuity in a narrow banded zone. Mathematically, strain localization is a bifurcation problem from a continuous deformation to a discontinuous one. Over the past decades, the major studies of shear banding are focused on dry or fully saturated geomaterials \cite{vardoulakis1995bifurcation}. However, in recent years, strain localization in unsaturated soils has been the subject of tremendous interests because of the importance of practical engineering applications (e.g., unsaturated soil slope failure and geo-energy production and storage) \cite{gens2010soil, morse2017evolution} ({see \cite{likos2019fundamental} for a recent review of this subject}). Unsaturated soils are three-phase porous media consisting of a deformable solid skeleton, pore water, and pore air \cite{alonso1990constitutive, fredlund1993soil, lu2004unsaturated}. Strain localization in unsaturated soils involves coupled multiphysics processes such as the coupled solid deformation and unsaturated fluid flow process (e.g., \cite{borja2013critical,song2014mathematical,song2014finite,song2018localized}). The numerical simulations via the mixed finite element method have demonstrated that material heterogeneities (e.g., density and the degree of saturation) have a first-order triggering role in shear banding in geomaterials \cite{song2017transient, song2017strain}. Multiphysics processes are known triggering agents in the failure of multiphase materials such as soils \citep{borja2013critical, song2014mathematical, lazari2015local}. The strain localization process usually presents multiscale and heterogeneous characteristics. Numerical methods such as the finite element method \cite{hughes2000finite, de2012nonlinear} are useful tools to study the onset of strain localization in geomaterials \citep{loret1991dynamic, lu2017multiscale}.  

The finite element method is formulated in terms of the weak or variational form of partial differential equations assuming continuous field variables such as displacement and fluid pressures \cite{lewis1998finite, zienkiewicz1999computational}. For this reason, the finite element method has difficulties in modeling localized failures/cracks because of the singularities across discontinuities \cite{stakgold2011green}. From the mathematical viewpoint, one underlying reason is the loss of ellipticity of the partial differential equation at the onset of extreme material failures (e.g., shear bands) \cite{deBorst1991simulation,de1993fundamental,de2012nonlinear}. Typical regularization strategies or techniques such as rate-dependent plastic models, higher-order gradient models, nonlocal models, the strong discontinuity approach, phase-field models, the extended finite element method have been formulated for single-phase solid materials to overcome the ill-posedness problem\cite{needleman1988material, de2020viscoplastic,de1992gradient, pijaudier1987nonlocal, simo1993analysis, borja2001strain, miehe2010phase, moes2002extended, belytschko2013nonlinear}. Recently these remedies for single-phase materials have been applied to simulate shear bands or fractures in multiphase materials \cite{deBorst1999coupled, ehlers2004deformation, oka2019computational, deBorst2006numerical, callari2010strong, mroginski2014discontinuous, mohammadnejad2013hydro, ehlers2017phase, roth2020fully}. For instance, R{\'e}thor{\'e} et al \cite{rethore2007two,rethore2007discrete} formulated a finite element model for simulating fluid flow in a fracturing unsaturated porous medium in which the discontinuity is described by exploiting the partition-of-unity property of finite element shape functions. Ehlers and Luo \cite{ehlers2017phase} proposed a phase-field approach embedded in the theory of porous media for modeling dynamic hydraulic fracturing. Oka et al \cite{oka2019computational} developed an elasto-viscoplastic model for strain localization analysis of unsaturated soils under dynamic conditions. Roth et al \cite{roth2020fully} presented a three-dimensional extended finite element model for application to complex flow in discontinuities.

Those aforementioned enhanced finite element methods (e.g., strong discontinuities) for modeling shear bands except the nonlocal models usually assume a zero thickness of shear bands. However, the thickness of shear bands in geomaterials can be in the order of several grain sizes as observed in the laboratory testing \cite{muhlhaus1987thickness, nemat1999micromechanics, rechenmacher2003digital, alshibli2008spatial} or be of the scale of meters in the field \cite{zhang1999interal, zhang2004particular, schrefler2006interaction}. Thus, in general, a length scale is required to realistically model the formation of shear bands at the continuum scale. Similarly, fluid flow in natural porous media implies that the flux at a given location may depend on a possibly large neighborhood that includes complex multiscale pathways of varying length and hydraulic conductivity. In general, there are both physical and numerical reasons for formulating a nonlocal framework for modeling variably saturated geomaterials. The physical reasons include: (i) heterogeneity of microstructures and its homogenization at a small scale on which the smoothed strain field cannot be considered as uniform \cite{bazant2002nonlocal}; (ii) grain or pore size and geometrical effects \cite{dormieux2006microporomechanics}; (iii) {length scales of modeling fluid flow in porous media at the continuum scale through the thermodynamically constrained averaging theory and upscaling methods \cite{gray2014introduction, cushman2013physics}; and (iv) solid-fluid interaction at the nanoscale considering long-range forces \cite{israelachvili2011intermolecular, song2019molecular}}. For example, de Borst et al \cite{deBorst1991simulation} presented a nonlocal mathematical model for strain localization using the Cosserrat continuum model of solid mechanics, wherein individual particles are enriched with three rotational degrees of freedom alongside the pre-existing three translational degrees of freedom. Zhang et al \cite{zhang1999interal} explored an internal length scale in multiphase materials such as variably saturated porous media and the relation of this length scale to the strain localization behavior under dynamic loading conditions. The nonlocal continuum concept has emerged as an effective means for regularizing the boundary value problems with strain softening, capturing the size effects, and avoiding spurious localization that gives rise to pathological mesh sensitivity in numerical computations \cite{de1993fundamental}. For instance, Lazari et al \cite{lazari2015local} employed nonlocal rate-dependent plasticity to regularize the numerical solution in the strain localization process of unsaturated geomaterials.

To model localized failure in unsaturated geomaterials, a promising novel approach is peridynamics. Peridynamics was proposed to unite the mathematical modeling of discontinuities, porous media, and discrete particles within the same mathematical framework \cite{silling2000reformulation, silling2007peridynamic}. Classical theory assumes that all bodies have a continuous, even distribution of mass and that all internal forces act at zero distance. In peridynamics, the governing equations are based on the concept that forces act on the material points over a finite distance. The classical theory falls apart at smaller length scales where long-range forces begin affecting material behavior. The field equations in peridynamics are integral-differential equations and are valid at discontinuities such as shear bands. Peridynamics has been applied to study the fluid flow and deformation problems in solid materials \cite{jabakhanji2015peridynamic, oterkus2017fully, ouchi2015fully, turner2013non, lai2015peridynamics, song2018peridynamics, menon2019coupled}. For instance, Turner \cite{turner2013non} proposed a peridynamic model accounting for fluid-structure interaction for modeling fractures. Jabakhanji et al \cite{jabakhanji2015peridynamic} proposed a bond-based peridynamic formulation for fluid flow in porous media. Katiyar et al \cite{katiyar2014peridynamic} formulated a state-based peridynamic model for fluid flow in geomaterials. Ouchi et al \cite{ouchi2015fully} presented a coupled state-based peridynamics formulation to simulate fluid-driven fractures in a poroelastic medium. In \cite{song2018peridynamics}, a peridynamic elastoplastic model was proposed for unsaturated geomaterials through the original correspondence principle \cite{silling2007peridynamic}. It is noted that the formulation in \cite{song2018peridynamics} is for the single-phase geomaterials because the unsaturated water flow equation is missing in the mathematical model.

In this article, we implement a computational periporomechanics model for unsaturated porous media in which the unsaturated fluid flow equation is coupled with the force balance equation. The mathematical model is formulated based on the peridynamic state concept \cite{silling2007peridynamic} and the peridynamic effective force concept for unsaturated porous media \cite{song2020peridynamic_unsaturated}.  The coupled field equations are integral-differential equations (integration in space and differentiation in time) without assuming the continuity of solid skeleton displacement and pore fluid pressures. The intrinsic multiphase length scale in the field equations can be utilized to capture the length scale characteristics associated with the multiphysics behavior of unsaturated porous media (e.g., thickness of deformation band and pore pressure band). Here the intrinsic multiphase length scale means that a length scale exists in both the solid skeleton deformation and fluid flow field equations. For numerical implementation, the governing equations are discretized in space by a multiphase Lagrangian meshfree method and in time by a fully implicit method \cite{hughes2000finite}. Due to its nonlocal formulation, the periporomechanics model is computationally more demanding than numerical methods based on classical local theory. To circumvent this problem, the implemented numerical model is enhanced for high-performance computing through open MPI (message passing interface) \cite{gropp1996high, gabriel2004open}. To validate the numerical implementation, numerical simulations of Terzaghi one-dimensional consolidation problem \cite{terzaghi1951theoretical}, drainage of a soil column \cite{liakopoulos1964transient}, and Mandel's slab problem \cite{mandel1953consolidation} are presented. Numerical simulations of localized failure in unsaturated porous media are performed through the implemented computational periporomechanics model. It is demonstrated by the detailed numerical results in Section \ref{numerical_examples} that the computational periporomechanics model is robust in simulating the formation of deformation and fluid flow bands in unsaturated porous media.

The remainder of the article is structured as follows. Section 2 introduces the unsaturated periporomechanics under quasi-static condition assuming a passive atmospheric pressure. Section 3 details the numerical implementation via an implicit multiphase Lagrangian meshfree method. Section 4 deals with validation of the numerical implementation. Section 5 presents numerical simulations of localized failure in unsaturated soils via the implemented computational periporomechanics model, followed by a summary in Section 6. Sign convention in continuum mechanics is adopted, i.e., the tensile stress in solid is positive, and pore fluid pressure in compression is positive.

\section{Mathematical model}
This section introduces the unsaturated periporomechanics model under quasi-static condition assuming passive air pressure. For the sake of notation and nomenclature, we briefly introduce the effective stress concept of classical unsaturated poromechanics (e.g., \cite{fredlund1993soil, lu2004unsaturated}). Unsaturated soil is a three-phase porous medium consisting of a solid skeleton, liquid water, and gaseous phases. The difference between pore air and pore water pressures is usually termed as matric suction. Here, we assume passive air pressure in unsaturated soils to be at equilibrium with external atmospheric pressure and thus can neglect its effect \cite{lewis1998finite}. In this case, matrix suction is merely negative pore water pressure. Here the generalized form of effective stress for the solid skeleton (e.g., \cite{nuth2008effective, lu2010closed}) is adopted.
\begin{equation}
\bar{\sigma}_{ij} = \sigma_{ij} + \left[(1-\bar{S}_r)p_a + \bar{S}_r p_w\right]\delta_{ij},
\label{general_effective_stress}
\end{equation}
where $\bar{\sigma}_{ij}$ is the effective Cauchy stress tensor, $\sigma_{ij}$ is the total Cauchy stress tensor, $p_a$ and $p_w$ are the pore air and pore water pressures respectively, $\delta_{ij}$ is the Kronecker delta, and $\bar{S}_r$ is the effective degree of saturation \cite{lu2010closed} defined as
\begin{equation}
\bar{S}_r = \frac{S_r - S_1}{S_2 - S_1},
\end{equation}
where $S_r$ is the degree of saturation, $S_1$ and $S_2$ is the residual and maximum degree of saturation, respectively.
Under the assumption of passive pore air pressure, equation \eqref{general_effective_stress} can be written as
\begin{equation}
\bar{\sigma}_{ij} = \sigma_{ij} + \bar{S}_r p_w\delta_{ij}.
\label{effective_stress}
\end{equation}



\subsection{Coupled unsaturated  periporomechanics model} \label{coupled_periporomechanics}

Peridynamics theory assumes a material is composed of material points, each interacting with its neighboring points within a nonlocal region \cite{silling2000reformulation, silling2007peridynamic}. Each interaction pair of a material point with its neighboring material point in the horizon is referred to as a ``bond". Peridynamic states are mathematical objects that can apply to continuous or discontinuous functions \cite{silling2007peridynamic}. Peridynamic states depend upon position and time and operate on a vector connecting any two material points. Depending on whether the value of this operation is a scalar or vector, a peridynamic state is called a scalar-state or a vector-state. For instance, the deformation state is a function that maps any mechanical bond onto its image under the deformation. A force state contains the forces within bonds of all lengths and orientations. In general, a peridynamic state can be viewed as a data set to extract information about the mechanical or physical state of material points. 

\begin{figure}[h!]
\centering
    \includegraphics[width=0.5\textwidth]{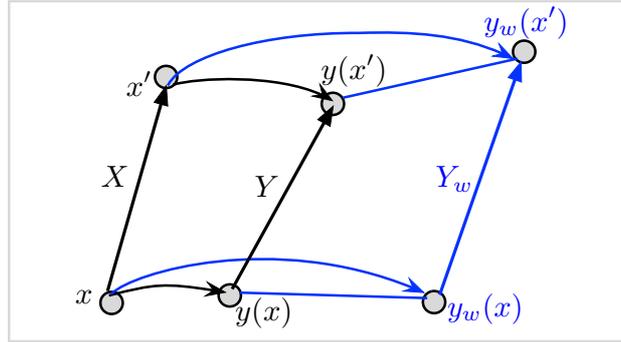}%
  \caption{Schematic of kinematics of the coupled periporomechanics assuming atmospheric air pressure.}
  \label{kinematics_peri_poromechanics}
\end{figure}

Figure \ref{kinematics_peri_poromechanics} sketches the kinematics of unsaturated periporomechanics in which pore air material points are not shown under the assumption of passive air pressure. The relative position of the solid material points $\vec{x}$ and $\vec{x}'$ in the reference configuration is given by the reference position vector $\vec{X}$ denoted as the mechanical bond vector $\vec{\xi}$ 
\begin{equation}
\vec{X}=\vec{x}'-\vec{x}=\vec{\xi},~~~\text{and}~\vec{\xi}'=\vec{x} - \vec{x}'.
\label{reference_state}
\end{equation}

The relative position of the same solid material points in the deformed configuration is given by the deformation position vector state $\vec{Y}$
\begin{equation}
\vec{Y} = \vec{y}'- \vec{y} ,
\label{deformation_state}
\end{equation}
where $\vec{y}'$ and $\vec{y}$ are the position vectors of the solid material points at positions $\vec{x}'$ and $\vec{x}$, respectively, in the deformed configuration. Let $\vec{u}=\vec{u}(\vec{x})$ and $\vec{u}'=\vec{u}(\vec{x}')$ are the displacements of material points at $\vec{x}$ and $\vec{x}'$, respectively. We have
\begin{equation}
\vec{y}=\vec{x}+\vec{u}, ~~~~\text{and} ~~~~\vec{y}'=\vec{x}'+\vec{u}'.
\end{equation}
The displacement vector states for bonds $\vec{\xi}$ and $\vec{\xi}'$ are defined as
\begin{equation}
\vec{U}=\vec{u}'- \vec{u}, ~~~~\text{and} ~~~~\vec{U}'=\vec{u} - \vec{u}'.
\label{displacement_state}
\end{equation}
Substituting the displacement vector state and the reference vector state at material point $\vec{x}$ into Equation \eqref{deformation_state} generates
\begin{equation}
\vec{Y}=\vec{X} + \vec{U}.
\end{equation}

\begin{figure}[h!]%
\centering
    \includegraphics[width=0.5\textwidth]{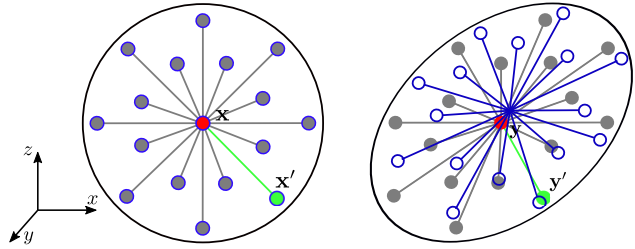}%
  \caption{Schematics of a solid material point (in gray and red)  and a pore water material point (in blue) and their corresponding neighboring material points within the horizon: Initial state (left) and deformed state (right).}
  \label{coupled_initial_deformed_states}
\end{figure}

In the unsaturated periporomechanics model, the Lagrangian coordinate system is adopted for the solid material points and the relative Eulerian coordinate system moving with respect to the solid phase is employed for the pore fluid material points (see Figure \ref{coupled_initial_deformed_states}). Thus, the relative position of the same pore water material points in the deformed configuration is given by the deformation position vector state $\vec{Y}_{w}$ 
\begin{equation}
\vec{Y}_{w} = \vec{y}'_{w} - \vec{y}_{w},
\end{equation}
where $\vec{y}'_{w}$ and $\vec{y}_{w}$ are the position vectors of the pore water material points at positions $\vec{x}'$ and $\vec{x}$, respectively, in the deformed configuration. The relative deformed position vector state of the pore fluid material points with respect to the solid material points is 
\begin{equation}
\widetilde{\vec{Y}}_{w} =\vec{Y}_{w} - \vec{Y}.
\end{equation}
Furthermore, we define the pore water pressure states at material points $\vec{x}$ and $\vec{x}'$ operating on the bonds, $\vec{\xi}$ and $\vec{\xi}'$, respectively, in the current configuration as
\begin{equation}
\Phi_w = p'_w - p_w, ~~\text{and}~~ \Phi'_w = p_w - p'_w, 
\end{equation}
where $p'_w$ and $p_w$ are pore water pressures at material points $\vec{x}'$ and $\vec{x}$, respectively.

Assuming a passive pore air pressure the effective force state $\overline{\vec{T}}$ can be defined as
\begin{equation}
\overline{\vec{T}}=\vec{T} - \bar{S}_r\vec{T}_{w},
\label{effective_force_state}
\end{equation}
where $\vec{T}$ is the total force state and $\vec{T}_w$ is the force state of pore water, which is related to pore water pressure. It has been demonstrated that the effective force state is an energy conjugate to the deformation state of the solid skeleton \cite{song2020peridynamic_unsaturated}. It follows from \eqref{effective_force_state} that the peridynamic linear momentum balance of unsaturated porous media under the quasi-static condition at $\vec{x}$ can be written as
\begin{equation}
\int_{\mathcal{H}_{\vec{x}}}\left[\left(\overline{\vec{T}} + \bar{S}_r\vec{T}_{w}\right) - \left(\overline{\vec{T}}' + \bar{S}'_r\vec{T}'_{w}\right)\right]\text{d}V'+\left[\rho_s(1-\phi)+\rho_w\phi \bar{S}_r\right]\vec{g}=\vec{0},
\label{peridynamic_momentum}
\end{equation}
where $\overline{\vec{T}} + \bar{S}_r\vec{T}_{w}$ and $\overline{\vec{T}}' + \bar{S}'_r\vec{T}'_{w}$ are the total force states at material points $\vec{x}$ and $\vec{x}'$, respectively, $\mathcal{H}_{\vec{x}}$ is the neighborhood of $\vec{x}$ with a characteristic length $\delta$, $\rho_s$ is the intrinsic density of the soil solid, $\rho_w$ is the intrinsic density of pore water, $\phi$ is the porosity of the mixture, and $\vec{g}$ is the gravity acceleration vector.

Let $\dot{\mathcal{V}}_s$ and $\dot{\mathcal{V}}'_s$ be the rates of volume change scalar state of the solid skeleton at material points $\vec{x}$ and $\vec{x}'$, respectively. Let $\widetilde{Q}_w$ and $\widetilde{Q}'_w$ be the mass flow scalar states of pore water relative to the solid skeleton at material points $\vec{x}$ and $\vec{x}'$, respectively. Assuming that the soil solid is incompressible, the peridynamic mass balance equation of the mixture can be written as
\begin{equation}
\phi\dfrac{\text{d}\bar{S}_r}{\text{d}t}+\dfrac{\phi \bar{S}_r}{K_w}\dfrac{\text{d}p_w}{\text{d}t}+\bar{S}_r\int_{\mathcal{H}_{\vec{x}}}\left(\dot{\mathcal{V}}_s - \dot{\mathcal{V}}'_s\right)\text{d}V'+\dfrac{1}{\rho_w}\int_{\mathcal{H}_{\vec{x}}}\left(\widetilde{Q}_w - \widetilde{Q}'_w \right)\text{d}V'=0,
\label{peridynamic_mass_balance}
\end{equation}
where $K_w$ is the bulk modulus of pore water.

\subsection{Constitutive and physical models }

In the proposed periporomechanics model, the constitutive and physical models for the solid deformation and unsaturated fluid flow are cast through the multiphase constitutive correspondence principle \cite{song2020peridynamic_unsaturated}. The nonlocal deformation gradient can be determined from the deformation state $\vec{Y}$ as
\begin{equation}
\widetilde{\vec{F}} = (\underline{\vec{Y}}*\underline{\vec{X}}){\vec{\mathcal{K}}}^{-1} = \left(\int_{\mathcal{H}_{\vec{x}}} \underline{\omega}\langle\vec{\xi}\rangle{\vec{Y}}\langle\vec{\xi}\rangle\otimes{\vec{X}}\langle\vec{\xi}\rangle\,dV_{x'}\right) \vec{\mathcal{K}}^{-1},
\end{equation}
where $\underline{\omega}\langle\vec{\xi}\rangle$ is a weighting function and  $\vec{\mathcal{K}}$ is the peridynamic shape tensor \cite{silling2007peridynamic} that is defined as 
\begin{equation}
\vec{\mathcal{K}} = \int_{\mathcal{H}_{\vec{x}}}\underline{\omega}\langle\boldsymbol\xi\rangle\vec{\xi}\otimes\vec{\xi}\,dV'.
\end{equation} 
The symbol `$*$' is the peridynamic equivalent of the tensor-product and represents a summation of the product of the two quantities over the entire horizon. The nonlocal deformation gradient represents a relative motion between a material point and its neighbors in the horizon.

Next, through the multiphase correspondence principle that the effective force state at material point $\vec{x}$ can be written as
\begin{equation}
\overline{\vec{T}}=\underline{\omega}\left(\vec{\sigma}+\bar{S}_rp_w\vec{1}\right)J\widetilde{\vec{F}}^{-T}\vec{\mathcal{K}}^{-1}\vec{\xi},
\label{effective_force_state_correspondence}
\end{equation}
where the weighting function is expressed as $\underline{\omega}$ for brevity, $\vec{1}$ is the second-order identity tensor, and $J$ is Jacobian (determinant) of $\widetilde{\vec{F}}$.
From \eqref{effective_force_state_correspondence} the total force state of the mixture reads
\begin{equation}
\vec{T}=\underline{\omega}\left(\overline{\vec{\sigma}}-\bar{S}_r p_w\vec{1}\right)J\widetilde{\vec{F}}^{-T}\vec{\mathcal{K}}^{-1}\vec{\xi}.
\label{total_force_state_correspondence}
\end{equation}
Given \eqref{total_force_state_correspondence}, the linear momentum balance equation of unsaturated porous media can be rewritten as
\begin{align}
&\int_{\mathcal{H}_{\vec{x}}}\left[\left(\underline{\omega}\left(\overline{\vec{\sigma}}-\bar{S}_rp_w\vec{1}\right)J\widetilde{\vec{F}}^{-T}\vec{\mathcal{K}}^{-1}\vec{\xi}\right) - \left(\underline{\omega}\left(\overline{\vec{\sigma}}'-\bar{S}'_rp'_w\vec{1}\right)J\widetilde{\vec{F}}^{-T}\vec{\mathcal{K}}^{-1}\vec{\xi}'\right)\right]\text{d}V'\nonumber\\
&+\left[\rho_s(1-\phi) + \rho_w\phi S_r\right]\vec{g}=\vec{0},
\label{linear_momentum_correspondence}
\end{align}
where the second term in the integrand is evaluated at material point $\vec{x}'$.  

The effective stress in equation \eqref{linear_momentum_correspondence} can be determined by a classical elastoplastic constitutive model for unsaturated soils through the multiphase constitutive correspondence principle. In unsaturated poromechanics, the degree of saturation can be determined through the soil water retention curve \cite{cao2018soil} given matric suction (i.e., negative pore water pressure). Without loss of generality, the soil water retention curve is represented by van Genuchten equation \cite{van1980closed}.
\begin{equation}
S_r = S_1 + (S_2-S_1)\left[1+\left(\frac{-p_w}{s_a}\right)^n\right]^{-m},
\label{SWCC}
\end{equation}
where $s_a$ is a scaling factor, and $m$ and $n$ are fitting parameters for the experimental data. 

For self-completeness, the key components of a classical elastoplastic constitutive model are summarized as follows. The total strain in the rate form reads,
\begin{equation}
\dot{\varepsilon}_{ij} = \dot{\varepsilon}_{ij}^e + \dot{\varepsilon}_{ij}^p,
\end{equation}
where $\varepsilon_{ij}$ is the total strain tensor, $\varepsilon_{ij}^e$ is the elastic strain tensor, $\varepsilon_{ij}^p$ is the plastic strain tensor, and $i, j= 1,2,3$. The effective stress is determined by the elastic model as
\begin{equation}
\dot{\bar{\sigma}}_{ij}=C_{ijkl}^e \dot{\varepsilon}_{kl},
\end{equation}
where $C_{ijkl}$ is the elastic stiffness tensor and $k,l=1,2,3$. 
\begin{equation}
C_{ijkl}^e=(K-2/3\mu)\delta_{ij}\delta_{kl} + \mu\left(\delta_{ik}\delta_{jl}+\delta_{il}\delta_{jk}\right),
\end{equation}
where $K$ and $\mu$ are the bulk and shear moduli of the solid skeleton, respectively. 
The plastic strain is determined by the plastic model. The yield surface is a function of the effective stress and the apparent pre-consolidation pressure $\bar{p}_c$ and is expressed as 
\begin{equation}
f(\bar{p},q,\overline{p}_c) =  \bar{p}(\bar{p} - \bar{p}_c) +  \frac{q^2}{M^2},
\end{equation}
where $\bar{p}$ is the mean effective stress, $\bar{p}_c$ is the apparent preconsolidation pressure, $q$ is the deviator stress, and $M$ is the slope of the critical state line. Here $\bar{p}$ and $q$ are determined by the following equations
\begin{equation}
\bar{p}  =\frac{1}{3}\bar{\sigma}_{ij}\delta_{ij} ~~\text{and}~~ q = \sqrt{\frac{3}{2}}\sqrt{s_{ij}s_{ij}}.
\end{equation}
The size of the yield function is dependent on the apparent preconsolidation pressure, which is a function of the plastic volumetric strain $\epsilon_v^p$, matric suction, and the degree of saturation \cite{borja2004cam, borja2013critical}. The plastic strain rate is determined by 
\begin{equation}
\dot{\varepsilon}_{ij}^p = \lambda\frac{\partial f}{\partial \bar{\sigma}_{ij}} = \lambda\left[ \frac{1}{3}(2\bar{p}- \bar{p}_c)\delta_{ij} + \frac{2q}{M^2}\sqrt{\frac{3}{2}}\hat{n}_{ij}\right]
\end{equation}
where $\lambda$ is the so-called plastic multiplier, and $\hat{n}_{ij} = s_{ij}/\sqrt{s_{kl}s_{kl}}$. 

Similarly, the non-local gradient of pore water pressure at the material point $\vec{x}$ in the current configuration can be determined from the pore water pressure state as
\begin{equation}
\widetilde{\vec{\nabla}\Phi}_w=\left(\int_{\mathcal{H}_{\vec{x}}} \underline{\omega}\langle\vec{\xi}\rangle{\Phi_w}\langle\vec{\xi}\rangle\vec{X}\langle\vec{\xi}\rangle\,dV'\right) \vec{\mathcal{K}}^{-1}.
\end{equation}
Given the nonlocal gradient of pore water pressure, the volume flux of the pore water $\vec{q}_w$ and $\vec{q}'_w$ at material point $\vec{x}$ and point $\vec{x}'$ can be determined by the generalized Darcy's law for unsaturated fluid flow \cite{bear2013dynamics} as
\begin{align}
{\vec{q}}_w=-k^r_w k\vec{1}\widetilde{\vec{\nabla}\Phi}_w~\text{and} ~{\vec{q}}'_w=-k^r_w k \vec{1}\widetilde{\vec{\nabla}\Phi}'_w,
\end{align}
where $k$ is the intrinsic hydraulic conductivity for pore water, $k^r_w$ is the relative permeability of pore water under unsaturated condition, and $\widetilde{\vec{\nabla}\Phi}'_w$ is determined at material point $\vec{x}'$. Given the soil water retention curve, the relative permeability of the pore water phase \cite{coussy2004poromechanics} is expressed as
\begin{equation}
k^r_w=\sqrt{S_r}\left[1-\left(1-S_r^{1/m}\right)^m\right]^2,
\end{equation}
where $m$ is the same parameter in Equation \eqref{SWCC}. Based on the multiphase correspondence principle for porous media \cite{song2020peridynamic_unsaturated}, the mass flow state for the pore water at material points $\vec{x}$ and $\vec{x}'$ can be defined as
\begin{align}
\widetilde{Q}_w = \underline{\omega}\rho_w J \vec{q}_w\vec{\mathcal{K}}^{-1}\vec{\xi},\label{flow_state_at_x}~\text{and}~
\widetilde{Q}'_w=\underline{\omega}\rho_w J \vec{q}'_w\vec{\mathcal{K}}^{-1}\vec{\xi}'.
\end{align}
It follows that the peridynamic mass balance of the fluid phase under unsaturated condition can be written as
\begin{equation}
\phi\dfrac{\text{d}S_r}{\text{d}t}+\dfrac{\phi S_r}{K_w}\dfrac{\text{d}p_w}{\text{d}t}+S_r\int_{\mathcal{H}_{\vec{x}}}\left(\underline{\omega}\dot{\vec{y}}\vec{\mathcal{K}}^{-1}\vec{\xi}-\underline{\omega}\dot{\vec{y}}'\vec{\mathcal{K}}^{-1}\vec{\xi}'\right)\text{d}V'+\int_{\mathcal{H}_{\vec{x}}}\left(\underline{\omega} J\vec{q}_w\vec{\mathcal{K}}^{-1}\vec{\xi}-\underline{\omega}J\vec{q}'_w\vec{\mathcal{K}}^{-1}\vec{\xi}' \right)\text{d}V'=0.
\label{mass_balance_correpondence_state_form}
\end{equation}
Th first two rate terms in equation \eqref{mass_balance_correpondence_state_form} are the same as in classical poromechanics and only the divergence terms are replaced by peridynamic formulations in line with a nonlocal vector calculus \cite{silling2007peridynamic, gunzburger2010nonlocal}. All the input material parameters in the constitutive and physical models for the solid skeleton and fluid phases are measurable variables.

\section{Numerical implementation}

This section concerns the numerical implementation of the aforementioned unsaturated periporomechanics model. It is assumed in the mathematical formulation of periporomechanics that a body of unsaturated porous media consists of equal-sized mixed material points. In line with this hypothesis, the integral-differential equations of periporomechanics can be discretized by a Lagrangian meshless method in spatial domain \cite{silling2006emu} and a fully implicit method in temporal domain \cite{hughes2000finite}. In the proposed method, the problem domain is discretized into a finite number of equal-sized mixed material points (i.e., pore fluid and solid skeleton). Each material point has two kinds of degrees of freedom, i.e., displacement and pore fluid pressure. Thus, the meshless method is different from the one in \cite{silling2005meshfree} or \cite{song2020peridynamic_unsaturated}, which is for modeling single-phase solid materials. We refer to the literature for an in-depth report on meshfree methods \cite{li2002meshfree, liu2005introduction}. The message passing interface is utilized to enhance the implemented computational periporomechanics model for high-performance computing \cite{gropp1996high, gabriel2004open}.

\subsection{Implicit time integration}
An implicit time integration for is employed for the unsaturated periporomechanics framework. The fundamental unknowns at each mixed material point are displacement vector $\vec{u}$ and pore water pressure $p_w$. Given the values of displacement vector and pore water pressure at time step $n$, the implicit time integration \cite{hughes2000finite, de2012nonlinear} is used to determine the displacement vector and pore water pressure at time step $n+1$. At time step $n+1$, we can write out the momentum balance as
\begin{equation}
\vec{\mathcal{G}}_{n+1} = \int_{\mathcal{H}_{\vec{x}}}\left[\left(\overline{\vec{T}}-\bar{S}_r\vec{T}_w\right) - \left(\overline{\vec{T}}'-\bar{S}'_r\vec{T}'_w\right)\right]\text{d}V'+\left[\rho_s(1-\phi)+\rho_w\phi \bar{S}_r\right]\vec{g}=\vec{0},
\label{momentum_balance_peri_n_plus_1}
\end{equation}
where the variables at mixed material points $\vec{x}$ and $\vec{x}'$ are evaluated at time step $n+1$, and the subscript $n+1$ is omitted for brevity. The effective force state at time step $n+1$ is determined by the effective stress through Equation \eqref{effective_force_state}. Given the nonlocal deformation gradient, the effective stress is computed through the elastoplastic constitutive model introduced in the previous section. In line with the global implicit time integration, this elastoplastic constitutive model is numerically implemented through an implicit algorithm (e.g., \cite{simo1998computational, nova2003constitutive, borja2013plasticity, song2014mathematical}). We refer to \cite{song2018modeling} for the implementation of the numerical algorithm.

The backward finite difference is adopted to approximate the rate terms of the mass balance equation in time domain. Then we have
\begin{align}
\widetilde{\mathcal{M}}_{n+1}=&\left(\dfrac{\phi \bar{S}_r}{K_w}-\phi\dfrac{\partial \bar{S}_r}{\partial p_w}\right)\dfrac{p_w-p_{w,n}}{\Delta t}+\bar{S}_r\int_{\mathcal{H}_{\vec{x}}}\left[\underline{\omega}\dfrac{\vec{y}-\vec{y}_{n}}{\Delta t}\vec{\mathcal{K}}^{-1}\vec{\xi}-\underline{\omega}\dfrac{\vec{y}'-\vec{y}'_{n}}{\Delta t}\vec{\mathcal{K}}^{-1}\vec{\xi}'\right]\text{d}V'\nonumber\\
&+\int_{\mathcal{H}_{\vec{x}}}\left[\underline{\omega}\vec{q}_w\vec{\mathcal{K}}^{-1}\vec{\xi}-\underline{\omega}\vec{q}'_w\vec{\mathcal{K}}^{-1}\vec{\xi}'\right]\text{d}V'=0,
\label{mass_balance_backward_time_difference}
\end{align}
where a subscript $n+1$ for variables at time step $n+1$ is omitted for brevity, and variables at time step $n$ are noted with a subscript $n$. Multiplying $\Delta t$ on both sides of Equation \eqref{mass_balance_backward_time_difference} generates
\begin{align}
\mathcal{M}_{n+1}=&\left(\dfrac{\phi S_r}{K_w}-\phi\dfrac{\partial S_r}{\partial p_w}\right)(p_w-p_{w,n})+S_r\int_{\mathcal{H}_{\vec{x}}}\left[\underline{\omega}(\vec{y}-\vec{y}_{n})\vec{\mathcal{K}}^{-1}\vec{\xi}-\underline{\omega}(\vec{y}'-\vec{y}'_{n})\vec{\mathcal{K}}^{-1}\vec{\xi}'\right]\text{d}V'\nonumber\\
&+{\Delta t}\int_{\mathcal{H}_{\vec{x}}}\left[\underline{\omega}\vec{q}_w\vec{\mathcal{K}}^{-1}\vec{\xi}-\underline{\omega}\vec{q}'_w\vec{\mathcal{K}}^{-1}\vec{\xi}'\right]\text{d}V'=0.
\label{mass_balance_peri_n_plus_1}
\end{align}

\subsection{Linearization}

\begin{figure}[h!]%
\centering
    \includegraphics[width=0.7\textwidth]{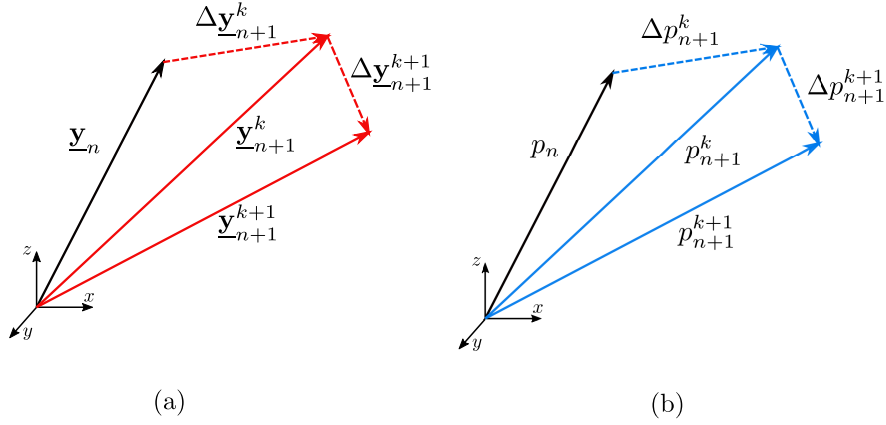}%
  \caption{Incremental kinematics for the solid phase and the pore water pressure of the fluid phase. Note: $k$ and $k+1$ are iteration numbers at time step $n+1$.}
  \label{incremental_kinematics}
\end{figure}

The tangent stiffness matrix is required when solving nonlinear equations by Newton's method. In general, the tangent stiffness matrix is a function of the chosen linearization point and must be updated every time step when the nodal degrees of freedom (displacement and pore water pressure) are updated. The tangent stiffness matrix relates an infinitesimal disturbance in the degree of freedom to change in the global force state and the global fluid flux state. For the state-based periporomechanics the tangent stiffness matrix contains non-zero entries for each node that interacts directly as well as nodes that share a neighbor material point, which makes the Jacobian inherently denser than that of a local poromechanics model. Figure \ref{incremental_kinematics} sketches the incremental kinematics for the solid phase and the pore fluid phase, respectively.

The linearization of the momentum balance equation with respect to the displacement vector states $\vec{U}$ and $\vec{U}'$, and pore water pressures $p_w$ and $p'_w$ at time step $n+1$ assuming an infinitesimal deformation is
\begin{align}
\delta\vec{\mathcal{G}}_{n+1} =& \int_{\mathcal{H}_{\vec{x}}}\left(\dfrac{\partial\overline{\vec{T}}}{\partial \vec{U}}\bullet\delta\vec{U}-\dfrac{\partial\overline{\vec{T}}'}{\partial\vec{U}'}\bullet\delta\vec{U}'\right)\text{d}V'\nonumber\\
&-\int_{\mathcal{H}_{\vec{x}}}\left[\left(\vec{T}_w\dfrac{\partial S_r}{\partial p_w}+S_r\dfrac{\partial\vec{T}_w}{\partial p_w}\right)\delta p_w-\left(\vec{T}'_w\dfrac{\partial S_r}{\partial p'_w}+S_r\dfrac{\partial\vec{T}'_w}{\partial p'_w}\right)\delta p'_w\right]\text{d}V'+\rho_w\phi\vec{g}\dfrac{\partial S_r}{\partial p_w}\delta p_w=\vec{0}.
\label{momentum_balance_peri_n_plus_1_linearization}
\end{align}

Similarly, the linearization of the mass balance equation at time step $n+1$ is

\begin{align}
\delta\mathcal{M}_{n+1}=&S_r\int_{\mathcal{H}_{\vec{x}}}\left[\underline{\omega}(\vec{\mathcal{K}}^{-1}\vec{\xi})\delta\vec{U}-\underline{\omega}(\vec{\mathcal{K}}^{-1}\vec{\xi}')\delta\vec{U}'\right]\text{d}V'+\left[\left(\dfrac{\phi S_r}{K_w}-\phi\dfrac{\partial S_r}{\partial p_w}\right) +\left(\dfrac{\phi}{K_w}\dfrac{\partial S_r}{\partial p}-\phi\dfrac{\partial^2 S_r}{\partial p^2_w}\right)\Delta p_w\right] \delta p_w\nonumber\\
&+\int_{\mathcal{H}_{\vec{x}}}\left[\underline{\omega}(\vec{\mathcal{K}}^{-1}\vec{\xi})\Delta\vec{y}-\underline{\omega}(\vec{\mathcal{K}}^{-1}\vec{\xi}')\Delta\vec{y}'\right]\text{d}V'\dfrac{\partial S_r}{\partial p_w}\delta p_w \nonumber\\
&+{\Delta t}\int_{\mathcal{H}_{\vec{x}}}\left[\underline{\omega}(\vec{\mathcal{K}}^{-1}\vec{\xi})\dfrac{\partial\vec{q}_w}{\partial \Phi_w}\delta \Phi_w-\underline{\omega}(\vec{\mathcal{K}}^{-1}\vec{\xi}')\dfrac{\partial\vec{q}'_w}{\partial \Phi'_w}\delta \Phi'_w\right)\text{d}V'=0.
\label{mass_balance_peri_n_plus_1_linearization}
\end{align} 

\subsection{Multiphase meshfree spatial discretization}

\begin{figure}[h!]%
\centering
    \includegraphics[width=0.7\textwidth]{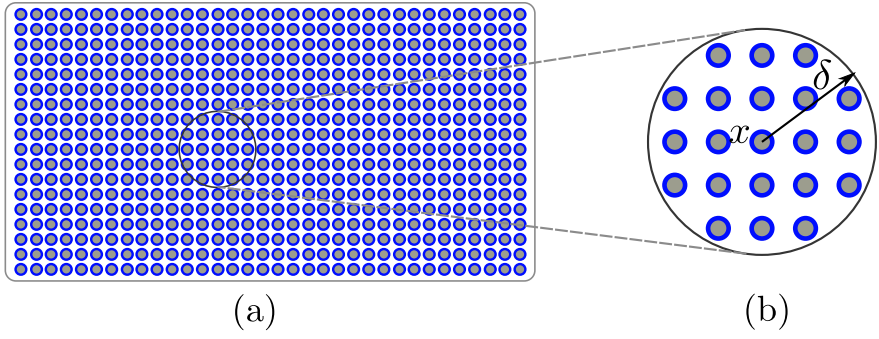}%
  \caption{(a) Schematic of the multiphase meshfree spatial discretization of an unsaturated porous material body, and (b) zoom-in view of material point $\vec{x}$ and its neighbors within the horizon $\delta$. Note: the gray dot represents the solid skeleton and the blue circle represents the pore water.}
  \label{meshfree_spatial_discretization}
\end{figure}

Figure \ref{meshfree_spatial_discretization} sketches the two-phase meshfree discretization of the problem domain.
For a periporomechanics material point $i$, there are totally $\mathscr{k}$ neighbor material points within its horizon. Here by the periporomechanics material point we mean that the material point has two states - the solid skeleton state and the pore fluid state. Therefore, each periporomechanics material point is endowed with two kinds of degree of freedom, the displacement of the solid skeleton and the pressure of the pore water. The momentum balance equation at material point $i$ becomes
\begin{equation}
\vec{\mathcal{G}}_{i} = \sum_{j=1}^{\mathscr{k}}\left[\left(\overline{\vec{T}}_i-\bar{S}_{r,i}\vec{T}_{w,i}\right) - \left(\overline{\vec{T}}'_j-\bar{S}'_{r,j}\vec{T}'_{w,j}\right)\right]V'_j+\left[\rho_s(1-\phi_i)+\rho_w\phi \bar{S}_{r,i}\right]\vec{g}=\vec{0},
\label{momentum_balance_peri_n_plus_1_spatial_discretization}
\end{equation}
where $V'_j$ is the volume of material point $j$ within the horizon of material point $i$.
Similarly, the mass balance equation at material point $i$ is
\begin{align}
\mathcal{M}_{i}=&\left(\dfrac{\phi \bar{S}_r}{K_w}-\phi\dfrac{\partial \bar{S}_r}{\partial p_w}\right)\Delta p_{w,i}+S_r\sum_{j=1}^{\mathscr{k}}\left[\underline{\omega}\Delta(\vec{y}_i)\vec{\mathcal{K}}^{-1}_i\vec{\xi}_i-\underline{\omega}(\Delta\vec{y}'_{j})\vec{\mathcal{K}}^{-1}_j\vec{\xi}'_j\right]V'_j\nonumber\\
&+{\Delta t}\sum_{j=1}^{\mathscr{k}}\left[\underline{\omega}\vec{q}_{w,i}\vec{\mathcal{K}}^{-1}_i\vec{\xi}_i-\underline{\omega}\vec{q}'_{w,j}\vec{\mathcal{K}}^{-1}_j\vec{\xi}'_j\right]V'_j=0.
\label{mass_balance_peri_n_plus_1_spatial_discretization}
\end{align}

At time step $n+1$, the linearization of the momentum balance is
\begin{align}
\delta\vec{\mathcal{G}}_{i} =&\sum_{j=1}^{\mathscr{k}}\left(\dfrac{\partial\overline{\vec{T}}_i}{\partial \vec{U}_i}\bullet\delta\vec{U}_i-\dfrac{\partial\overline{\vec{T}}'_j}{\partial\vec{U}'_j}\bullet\delta\vec{U}_j'\right)V'_j+\rho_w\phi_i\vec{g}\dfrac{\partial \bar{S}_{r,i}}{\partial p_{w,i}}\delta p_{w,i}\nonumber\\
&-\sum_{j=1}^{\mathscr{k}}\left[\left(\vec{T}_{w,i}\dfrac{\partial \bar{S}_{r,i}}{\partial p_i}+\bar{S}_{r,i}\dfrac{\partial\vec{T}_{w,i}}{\partial p_{w,i}}\right)\delta p_{w,i}-\left(\vec{T}'_{w,j}\dfrac{\partial \bar{S}_{r,j}}{\partial p'_{w,j}}+\bar{S}_{r,j}\dfrac{\partial\vec{T}'_{w,j}}{\partial p'_{w,j}}\right)\delta p'_{w,j}\right]V'_j=\vec{0}.
\label{momentum_balance_peri_n_plus_1_linearization_spatial_discretizatioin}
\end{align}

Similarly, the spatial discretization of the linearized mass balance equation is

\begin{align}
\delta\mathcal{M}_{i}=&\bar{S}_r\sum_{j=1}^{\mathscr{k}}\left[\underline{\omega}(\vec{\mathcal{K}}^{-1}_i\vec{\xi}_i)\delta\vec{y}_i-\underline{\omega}(\vec{\mathcal{K}}^{-1}_j\vec{\xi}'_j)\delta\vec{y}'_j\right]V'_j+\left[\left(\dfrac{\phi \bar{S}_r}{K_w}-\phi\dfrac{\partial \bar{S}_r}{\partial p_w}\right) +\left(\dfrac{\phi}{K_w}\dfrac{\partial \bar{S}_r}{\partial p_w}-\phi\dfrac{\partial^2 \bar{S}_r}{\partial p^2_w}\right)\Delta p_w\right] \delta p_{w,i}\nonumber\\
&+\sum_{i}^{\mathscr{k}}\left[\underline{\omega}(\vec{\mathcal{K}}^{-1}_i\vec{\xi}_i)\Delta\vec{y}_i-\underline{\omega}(\vec{\mathcal{K}}^{-1}_j\vec{\xi}'_j)\Delta\vec{y}'_j\right]V'_j\dfrac{\partial \bar{S}_r}{\partial p_w}\delta p_{w,i} \nonumber\\
&+{\Delta t}\sum_{j=1}^{\mathscr{k}}\left[\underline{\omega}(\vec{\mathcal{K}}^{-1}_i\vec{\xi}_i)\dfrac{\partial\vec{q}_{w,i}}{\partial \Phi_{w,i}}\bullet\delta \Phi_{w,i}-\underline{\omega}(\vec{\mathcal{K}}^{-1}_j\vec{\xi}'_j)\dfrac{\partial\vec{q}'_{w,j}}{\partial \Phi'_{w,j}}\bullet\delta \Phi'_{w,j}\right]V'_j=0.
\label{mass_balance_peri_n_plus_1_linearization_spatial_discretization}
\end{align}
The tangent matrix must be assembled by looping over all the individual particles and neighbors and account for the effect of perturbations in all degrees of freedom in every single bond on the force state and the fluid flux state at the point. The linearization must be considered for the periporomechanics material point itself as well as all its neighbors for each degree of freedom such as the displacement and the pore water pressure. Parallel computing capability of the implemented numerical model is realized through Open MPI, an open-source MPI (message passing interface) library \cite{silling2005meshfree, silling2006emu}. Here MPI is a communication protocol for programming parallel computers. We refer to \cite{gropp1996high, gabriel2004open} for the technical details regarding the implementation of Open MPI for high-performance computing.

\section{Validation of the numerical implementation}
In this section we validate the implemented computational periporomechanics model by conducting numerical simulations of classical poromechanics problems including Terzaghi's one-dimensional consolidation problem, the drainage of a soil column, and Mandel's slab problem \cite{terzaghi1951theoretical, liakopoulos1964transient, mandel1953consolidation}.

\subsection{Terzaghi's one-dimensional consolidation} \label{terzaghi_example}

\begin{figure}[h!]
\centering
       \includegraphics[width=0.12\textwidth, height=0.23\textheight]{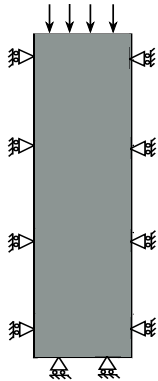}
  \caption{Schematic of the geometry and boundary conditions of the Terzaghi consolidation problem.}
  \label{1D_Setup}
\end{figure}

\begin{figure}[h!]
\centering
    \includegraphics[width=0.8\textwidth]{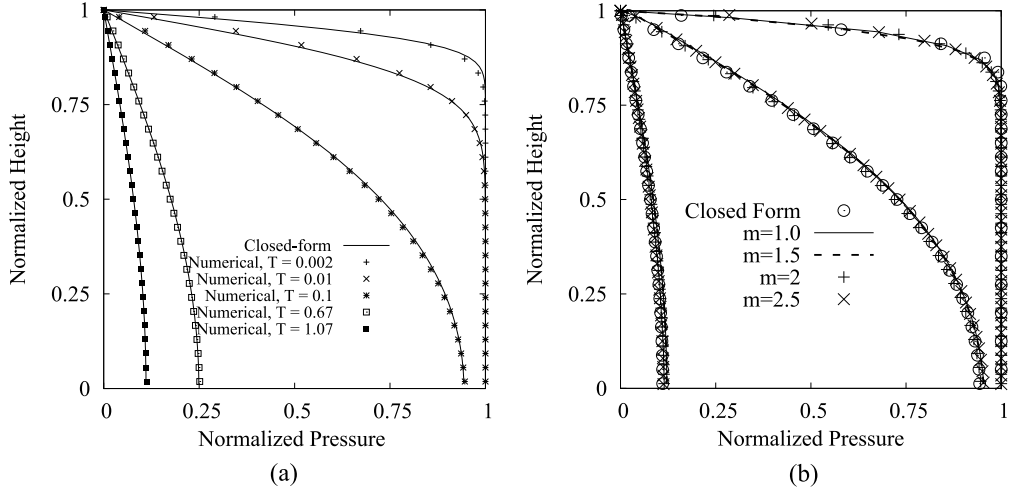}
  \caption{(a) Variations of the excess water pressure versus the height of the soil specimen from the numerical and closed-form solutions at different consolidation times. (b) Convergence of the periporomechanics solution to the closed-form solution with increasing nonlocal parameter. Note: T is a dimensionless time factor \cite{terzaghi1951theoretical}.}
  \label{1d_consolidation_results_comparison}
\end{figure}

The Terzaghi's one-dimensional consolidation problem is simulated, and the numerical result is compared against the closed-form solution. For the specimen, the height is 20 cm, the width is 5 cm, and the out-of-plane thickness is 1 cm.  The soil skeleton is assumed as a homogeneous elastic material with bulk modulus $K = 8.3 \times 10^4$ kPa and shear modulus $G$ = 1.8 $\times 10^4$ kPa. The fluid density $\rho_w$ is 1000 kg/m$^3$ and is assumed incompressible. Hydraulic conductivity of the soil skeleton $k$ is $1\times 10^{-8}$ m/s and the initial porosity is $\phi$ = 0.5. Figure \ref{1D_Setup} sketches the model set-up and the boundary conditions of the fully saturated soil specimen. Both lateral boundaries are constrained against horizontal deformation and are free to deform vertically. The bottom boundary is fixed in all three directions.  The problem domain is discretised into 400 material points with a characteristic dimension of $\Delta x$ = 0.5 cm. The horizon $\delta$ is assumed as $2\Delta x$.
The initial and boundary conditions for the fluid phase are as follows. At the initial state, all boundaries for the fluid phase are impervious. The top boundary is prescribed with a zero fluid pressure (i.e., the atmospheric pressure) during the consolidation process. A surcharge in the form of body force, equivalent to 50 kPa, is applied to the top boundary at t = 0 s and is held constant during the simulation. The total simulation time is 2000 s with an equal time step $\Delta t$ = 1 s. Figure \ref{1d_consolidation_results_comparison}(a) plots the excess pore water pressure from the numerical simulation against the closed-form solution \cite{terzaghi1951theoretical} at different consolidation stages.  Figure \ref{1d_consolidation_results_comparison}(b) shows the convergence of the numerical simulation to the analytical one with an increasing nonlocal parameter `m': the ratio of the horizon to the characteristic length of the material point (i.e., $\Delta x$). It is demonstrated from Figure \ref{1d_consolidation_results_comparison} that the numerical solution is in good agreement with the closed-form solution.

\subsection{Drainage of a soil column}

In this example, the experimental test of drainage of a sand column \cite{liakopoulos1964transient} is simulated numerically. The numerical results are compared with the experimental data. The sand column is 1 m in height and 0.25 m in width. The sample is discretized into uniform material points with a characteristic length of $\Delta x$ = 0.2 m. The horizon $\delta$ is assumed as $2\Delta x$. For the material parameters, Young's modulus is 1.3 MPa, Poisson's ratio is 0.4, soil density is 2000 kg/m$^3$, water density is 1000 kg/m$^3$, porosity is 0.2975, intrinsic permeability is 4.5 $\times$ 10$^{-10}$m$^2$.  From the data provided by Liakopoulos, the degree of saturation and the relative permeability of water are determined by the following equations \cite{lewis1998finite}.
\begin{align}
S_{r} & = 1 - (1.9722 \times 10^{-11})(-p_w)^{2.4279},\\
k_{r} &= 1-2.207 (1- S_{r})^{1.021}.
\end{align}

\begin{figure}[h!]
\centering
    \includegraphics[width=0.4\textwidth]{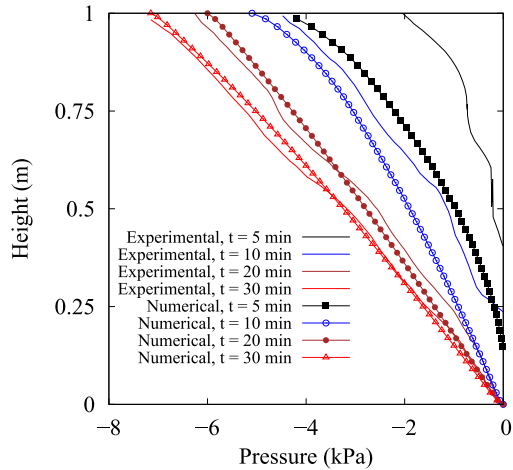}
  \caption{Comparison of the pore water pressures from the numerical simulation and experimental data \cite{liakopoulos1964transient} in the soil column at different times.}
  \label{drainage-comparison-numerical-experimental}
\end{figure}

\begin{figure}[h!]
\centering
    \includegraphics[width=0.4\textwidth]{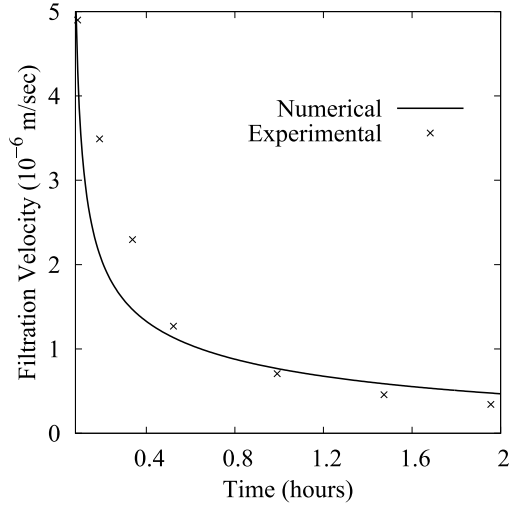}
  \caption{Comparison of the fluid flow velocities from the numerical simulation and experimental data \cite{liakopoulos1964transient} at the drainage boundary over time.}
  \label{filtration-comparison}
\end{figure}

For the solid phase, the lateral surfaces are constrained against lateral displacement but free to deform vertically, and the bottom is constrained for vertical displacement. For the fluid phase, both the upper and lower ends of the column are exposed to atmospheric pressure, and zero flux conditions are imposed on the lateral surfaces. The sand column is initially fully saturated with water. The water drainage is solely driven by the gravity load. The results are plotted in Figures \ref{drainage-comparison-numerical-experimental} and \ref{filtration-comparison}. In Figure \ref{drainage-comparison-numerical-experimental}, the pore water pressure from numerical simulations at different times is compared to the experimental results. The numerical results are in good agreement with the experimental data except for the early-stage results. In Figure \ref{filtration-comparison}, the fluid flow velocity at the drainage boundary predicted by the numerical simulation is compared against the experimental result. It is shown that the numerical result is close to the experimental data. It shall be noted that because the input parameters adopted from the literature are approximate ones the numerical results can not match the experimental data exactly.

\subsection{Mandel's slab problem}

\begin{figure}[h!]
\centering
    \includegraphics[width=0.4\textwidth]{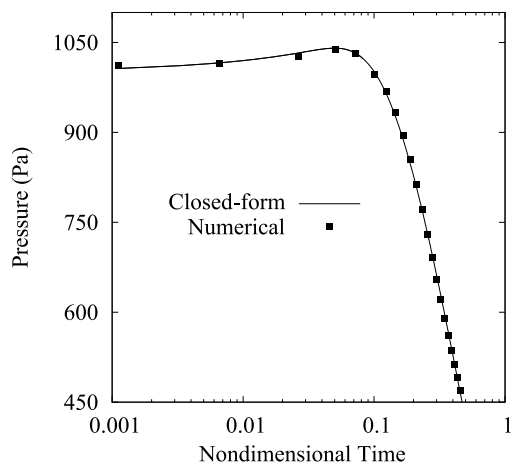}
  \caption{Comparison of the numerical and closed-form solutions of the water pressure at point A at the center of Mandel's slab. }
  \label{mandel-water-pressure-center}
\end{figure}

In this example, Mandel's slab problem is simulated to examine if the implemented computational periporomechanics model can capture the so-called Mandel-Cryer effect \cite{biot1941general, coussy2004poromechanics} under saturated condition. Here the Mandel-Cryer is a non-monotonic pore fluid pressure response that can only be captured by a fully coupled poromechanics solution. This phenomenon was first studied by Mandel with a saturated isotropic poroelastic specimen \cite{mandel1953consolidation} and later by Cryer with an isotropic poroelastic sphere under hydrostatic pressure \cite{cryer1963comparison}. Both studies showed that the induced pore fluid pressure at the center of the specimen is higher than the applied pressure during the initial loading stage. The physical reason behind this phenomenon is that fluid flow and pressure dissipation are retarded by the permeability and the distance to the drainage boundary. In contrast, the applied load and pore pressure response may be immediate and constant. Thus, when the pore fluid pressure near the drainage boundary starts to dissipate, the consequent rapid and localized deformation leads to an instantaneous load transfer to other parts of the porous body, which may lead to the development of additional excess pore pressure.
\begin{figure}[h!]
\centering
    \includegraphics[width=0.4\textwidth]{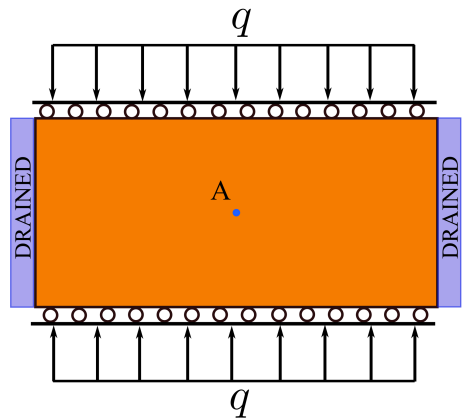}
  \caption{Sketch of the model setup, loading protocol and the boundary conditions. Note: the rollers are free to move horizontally and $q$ stands for the external load. }
  \label{mandel_setup}
\end{figure}%
For Mandel's slab problem, an infinitely long rectangular saturated soil specimen with free drainage conditions at both lateral boundaries is suddenly subjected to a constant vertical surcharge through a rigid and friction-free plate of the same width as the soil sample \cite{mandel1953consolidation} (see Figure \ref{mandel_setup}). The same material parameters for the Terzaghi's consolidation in \ref{terzaghi_example} are adopted. The surcharge is equal to the initial excess pore pressure in the specimen. The lateral boundaries are stress-free and can deform freely. The deformation of the specimen is under two-dimensional conditions. The dimension of the specimen is 10 m in width, 3 m in height, and 0.4 m in out-of-plane thickness. Figure \ref{mandel_setup} sketches the problem domain. The specimen is discretized into 1500 equal-sized material points with a characteristic length scale of $\Delta x$ = 0.2 m. The horizon is assumed as $\delta$ = 0.4 m. The initial excess pore fluid pressure is 1 kPa throughout the specimen. At t = 0, both lateral boundaries are forced to zero pressure, and a constant vertical surcharge equivalent to 1 kPa is applied to the top boundary. 
The numerical result is compared against the closed-form solution \cite{mandel1953consolidation,verruijt2013theory}. Figure \ref{mandel-water-pressure-center} plots the variation of the water pressure at the center of the Mandel's slab obtained by the numerical simulation and the closed-form solution. It is shown that the numerical result generated by the proposed computational periporomechanics model are in good agreement with the closed-form solution. Thus, it can be concluded that the fully coupled computational periporomechanics model can capture the Mandel-Cryer effect in the classical coupled poromechanics problem.

\section{Numerical examples} \label{numerical_examples}

Numerical examples of localized failure in unsaturated soils under different scenarios are presented in this section. Example 1 is a shear and compression test of a two-dimensional unsaturated soil sample. Example 2 concerns strain localization in a simple compression test under two-dimensional condition. Examples 1 and 2 study the influence of the horizon on shear bands. Example 2 also analyzes the impact of the intrinsic hydraulic conductivity and spatial discretization on the numerical results. Example 3 simulates the bearing capacity of a footing on unsaturated soils. For all examples, the same intrinsic length scales (i.e., horizon) are adopted for both the fluid and solid phases. For simplicity, the unit weighting function is used for all terms in the mathematical model.

\subsection{Example 1: Shear and compression of a two-dimensional specimen} \label{example_1}\label{example_1_base} 

\begin{figure}[h!]
\centering
    \includegraphics[width=0.60\textwidth]{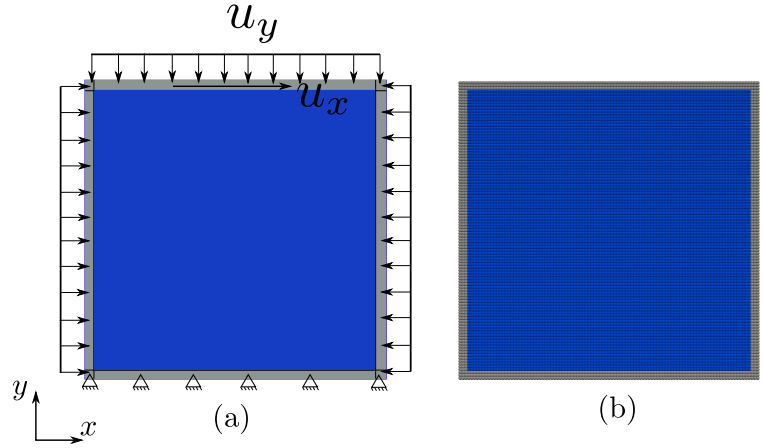}
  \caption{(a) Model setup, and (b) spatial discretization (the problem domain in blue and the boundary layers in gray).}
  \label{shear-setup}
\end{figure}

In this example, we simulate a two-dimensional unsaturated soil specimen under combined shear and compression loading. As shown in Figure \ref{shear-setup} (a), the soil specimen is 1 cm by 1 cm in the $x-y$ plane and the out-of-plane thickness is 0.025 cm. The sample is discretized into 12,800 material points with a characteristic length $\Delta x$ = 0.125 mm. Figure \ref{shear-setup}(b) sketches the spatial discretization. The horizon $\delta$ is 0.25 mm and the thickness of boundary layers are 0.25 mm. The material parameters for the solid skeleton are assumed as follows: bulk modulus is 38.9 MPa, shear modulus is 13 MPa, swelling index is 0.03, slope of critical state line is 1.0, compression index is 0.12, and specific volume at unit preconsolidation pressure is 2.76, and fitting parameters for preconsolidation pressure \cite{song2018peridynamics} are 0.185 and 1.42. The parameters for soil water characteristic curve are:  $S_1$ = 0.0,  $S_{2}$ = 1.0, $s_a$ = 20 kPa, $n$ = 2.0, $m$ = 0.5. The soil density is 2000 kg/$\mathrm{m}^3$, water density is 1000 kg/$\mathrm{m}^3$, initial porosity is 0.5, and hydraulic conductivity is 1$\times 10^{-8}$ m/s. The same material parameters are also adopted for the numerical examples presented in Section \ref{example_2}.

The bottom of the specimen is fixed. All fluid phase boundaries are no-flux conditions. For the initial state, a uniform isotropic stress state with a magnitude of -100 kPa is prescribed on the skeleton, and apparent pre-consolidation pressure is -300 kPa. In the sample, the initial suction is assumed 20 kPa, and the degree of saturation is 0.8. The constant confining pressure on the lateral boundaries is 84 kPa, which is imposed through an equivalent total force density \cite{song2018modeling}.  For the loading protocol, the soil specimen is subjected to a vertical compression load and an in-plane shear load on the top boundary. Both loads are prescribed by constant displacement rates of $1\times10^{-5}$ m/s in the $y$ direction and $1.4\times10^{-5}$ m/s in the $x$ direction, respectively (see Figure \ref{shear-setup}). {The confining pressure is applied on both lateral boundaries by a linear ramp function in 0.04 s. Then, the external loads are imposed on the top boundary with the prescribed constant rates}. The total simulation time is 24 s. 

\begin{figure}[h!]%
\centering
    \includegraphics[width=0.35\textwidth]{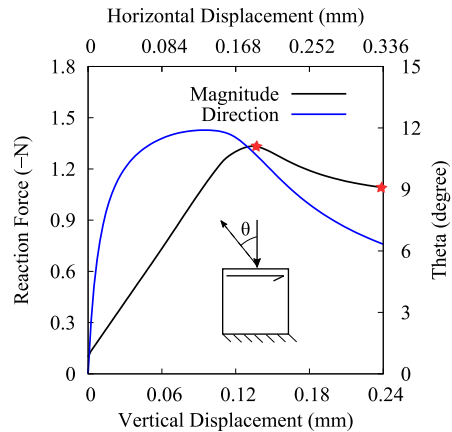}%
  \caption{Reaction force and its direction versus the vertical and horizontal displacements at the top boundary. Note the markers at displacements:  $u_y$ = 0.13 mm and  $u_y$ =  0.24 mm}
  \label{shear-force-plot}
\end{figure}%

\begin{figure}[h!]%
\centering
    \includegraphics[width=0.60\textwidth]{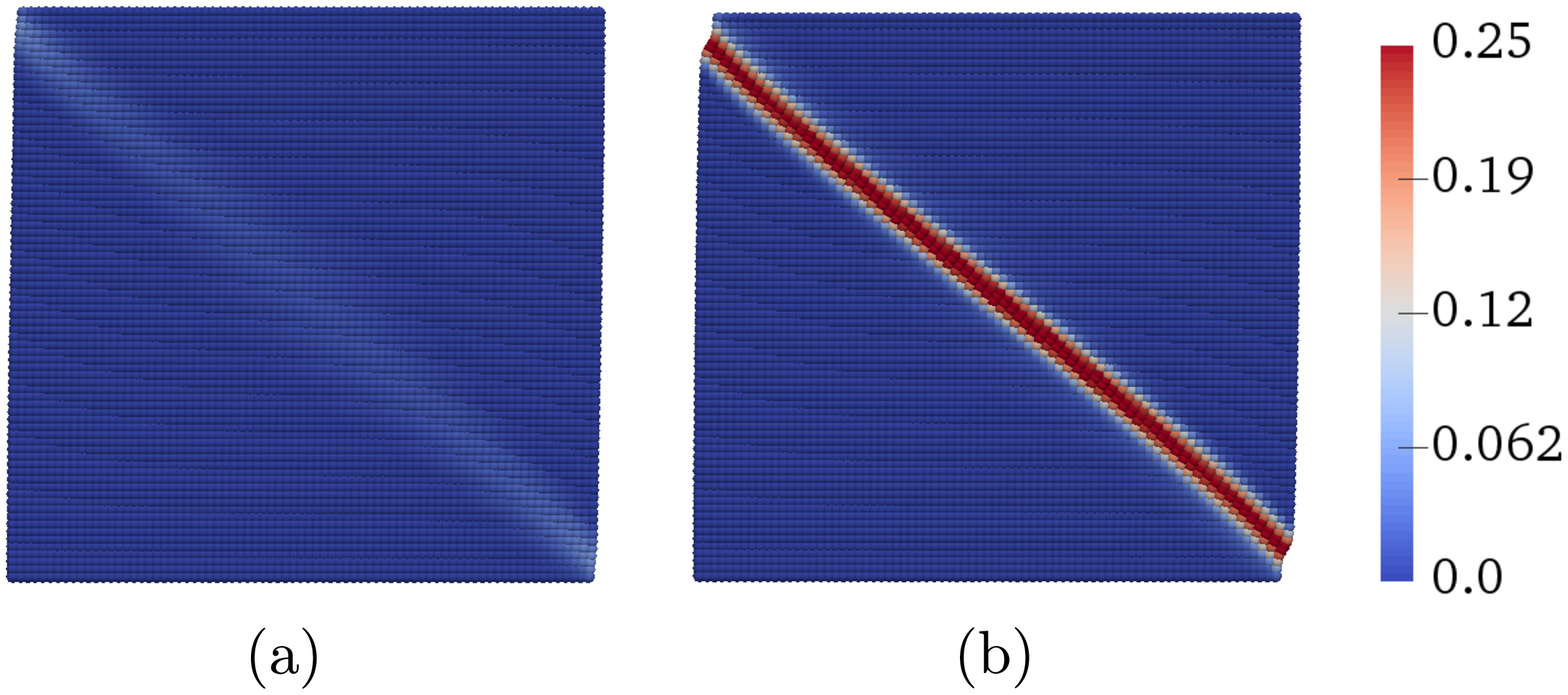}%
  \caption{Contours of equivalent plastic shear strain at displacements: (a) $u_y$ = 0.13 mm and (b) $u_y$ =  0.24 mm.}
  \label{shear-pds-Tcompare}
\end{figure}%

\begin{figure}[h!]%
\centering
    \includegraphics[width=0.60\textwidth]{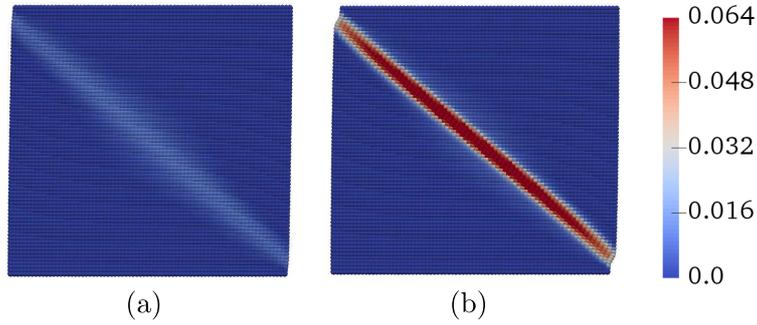}%
  \caption{Contours of plastic volumetric strain at displacements: (a) $u_y$ = 0.13 mm and (b) $u_y$ =  0.24 mm.}
  \label{shear-pdv-Tcompare}
\end{figure}%

\begin{figure}[h!]
\centering
    \includegraphics[width=0.60\textwidth]{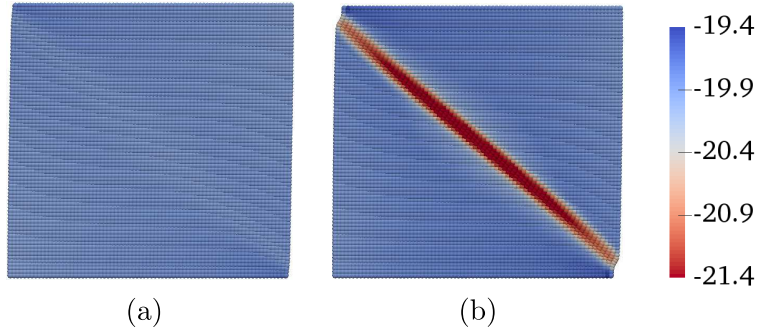}
  \caption{Contours of pore water pressure (kPa) at displacements: (a) $u_y$ = 0.13 mm and (b) $u_y$ =  0.24 mm}
  \label{shear-pressure-Tcompare}
\end{figure}

\begin{figure}[h!]%
\centering
    \includegraphics[width=0.60\textwidth]{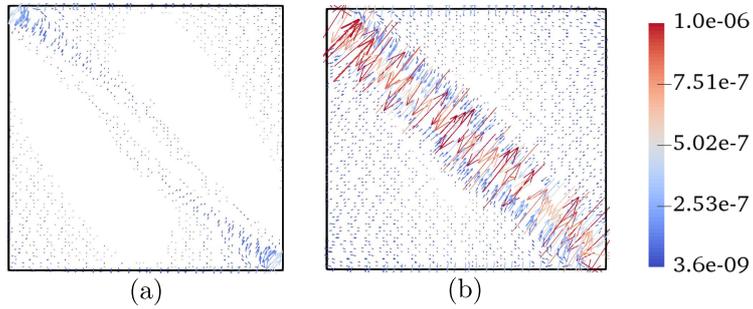}%
  \caption{Contours of fluid flux vector and the magnitude at displacements: (a) $u_y$ = 0.13 mm and (b) $u_y$ =  0.24 mm}
  \label{Example1-Flux-Base}
\end{figure}

Figure \ref{shear-force-plot} plots the reaction force versus the vertical and horizontal displacements at the top boundary, as well as {the direction of the reaction force and its eventual variation with increasing displacement.} The markers in Figure \ref{shear-force-plot} denote the vertical displacements: $u_y$ = 0.14 mm and $u_y$ = 0.24 mm, respectively. The results of plastic strains, pore water pressure  and fluid flux vector for the load steps are shown in Figures \ref{shear-pds-Tcompare}, \ref{shear-pdv-Tcompare}, \ref{shear-pressure-Tcompare} and \ref{Example1-Flux-Base}. Figure \ref{shear-pds-Tcompare} shows the contours of equivalent plastic shear strain (i.e., the second invariant of the plastic deviator strain tensor \cite{borja2013critical}). It is noted that the boundary layers are omitted here. Figure \ref{shear-pdv-Tcompare} plots the contours of plastic volume strain. The results in both figures demonstrate that a single shear band with a finite thickness has been initiated and developed in the specimen. The soil skeleton experiences positive plastic volumetric strain in the localized shear zone. Thus, it can be concluded that the shear band is a dilative one. Figure \ref{shear-pressure-Tcompare} plots the pore water pressure at the selected load steps. The results in Figure \ref{shear-pressure-Tcompare} show that a single banded zone has formed in the contour of pore water pressure. The banded zone of pore water pressure is consistent with the shear band of the solid skeleton shown in Figures \ref{shear-pds-Tcompare} and \ref{shear-pdv-Tcompare}. The decrease in the pore water pressure inside the shear band (equivalently, suction increases) can be in part due to the plastic dilation in the banded deformation zone. 

Figure \ref{Example1-Flux-Base} plots the snapshots of the fluid flux vector on the undeformed configuration at the selected load steps. The results demonstrate that the water flows into the center of the specimen at the inception of shear banding. After the shear band occurs, more water flows into the banded zone. This can be due to the fact that the pore water pressure in the banded zone is lower than the pore water pressure outside the shear band (refer to Figure \ref{shear-pressure-Tcompare}). As the plastic volumetric strain (see Figure \ref{shear-pdv-Tcompare}) increases in the banded zone, the pore water pressure decreases and the magnitude of the water flux rises because of the increased difference (or discontinuity) between the pore water pressures within and outside the shear band. It can be concluded from the example that the proposed periporomechanics model can capture the finite thickness of the shear band and the localized zone of the water pressure in the specimen.

\subsubsection{Influence of horizon}

\begin{figure}[h!]%
\centering
    \includegraphics[width=0.35\textwidth]{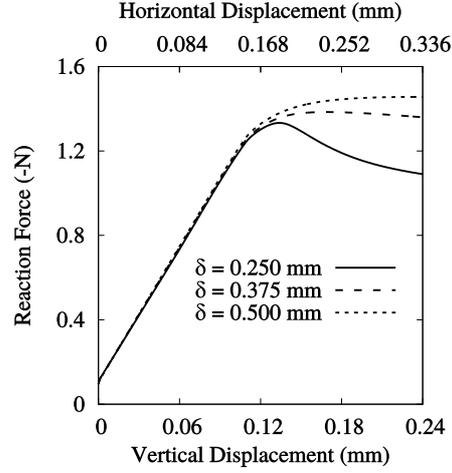}%
  \caption{Comparison of reaction force for the cases of $\delta$ = 0.25 mm, $\delta$ = 0.375 mm and $\delta$ = 0.5 mm. }
  \label{force-ShearPDS-Horzcompare}
\end{figure}

\begin{figure}[h!]%
\centering
    \includegraphics[width=0.60\textwidth]{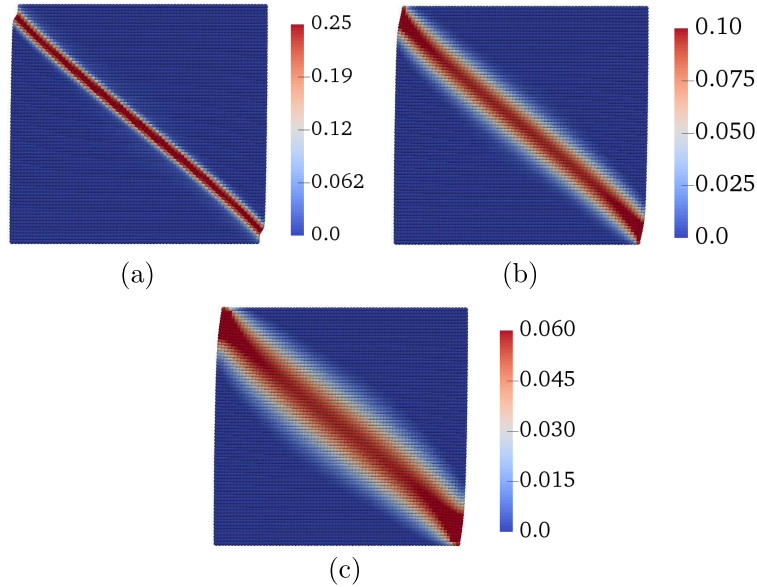}%
  \caption{Contours of equivalent plastic shear strain at the applied displacement $u_y$ = 0.24 mm for (a) $\delta$  = 0.25  mm, (b) $\delta$ = 0.375 mm,  and (c) $\delta$ = 0.5 mm.}
  \label{shear-PDS-Horzcompare}
\end{figure}

The horizon $\delta$  in the proposed periporomechanics is an intrinsic length scale for both the fluid and solid phases. The impact of the horizon on the formation of shear banding under shear and compression conditions is examined here. In so doing,  we re-run the simulation in Section \ref{example_1_base} with two more horizons, $\delta$ = 0.375 mm and $\delta$ =0.5 mm while keeping all other simulation parameters unchanged.  Figure \ref{force-ShearPDS-Horzcompare} compares the loading capacity curves (the reaction force versus displacement) generated by the simulations with the three horizons, respectively. The results show that the peak reaction forces (or applied loads) are different for the simulations with these three horizons. It is evident that the force responses are almost identical at the early loading stages. After the peak load, the loading curves diverge, which demonstrates the load capacity of the specimen in the post-localized regime is sensitive to the intrinsic length scale (horizon). Specifically, the simulation with a larger horizon (e.g., $\delta$ =0.5 mm) generates a relatively larger loading capacity in the post-localization regime. Figures \ref{shear-PDS-Horzcompare} and \ref{shear-pressure-Horzcompare} plot the contours of equivalent plastic shear strain and pore water pressure at the vertical displacement = 0.24 mm, respectively. The results obviously show that the internal length scale impacts both the formation of shear deformation band and pore pressure band. The simulation with $\delta$ = 0.25 mm (the smallest horizon in all three cases) produces the narrowest shear band and the maximum intensity of strain and pore water pressure within the shear band in all three cases. From the results in Figures \ref{force-ShearPDS-Horzcompare} - \ref{shear-pressure-Horzcompare}, we may conclude that both the deformation pattern, pore water pressure distribution, the load capacity of the specimen are dependent on the horizon size in the numerical model. For instance, under the same other conditions, the simulation with a larger horizon generates less intense plastic deformation in the localized zone than the simulation with a smaller horizon. Furthermore, the results show that the horizon may be correlated to the width of the shear band. It may imply that the proposed periporomechanics model can better capture the thickness of shear bands in unsaturated soils by adjusting the horizon. 

\begin{figure}[h!]%
\centering
    \includegraphics[width=0.60\textwidth]{Example1-Pressure-Horizons.eps}%
  \caption{Contours of pore water pressure (kPa) at the applied displacement $u_y$ = 0.24 mm for (a) $\delta$  = 0.25 mm ,(b) $\delta$ = 0.375 mm, and  (c) $\delta$ = 0.5 mm.}
  \label{shear-pressure-Horzcompare}
\end{figure}%

\subsection{Example 2: Compression of a two-dimensional specimen} \label{example_2}

In this example, we conduct a compression test of a two-dimensional unsaturated soil specimen with a constant confining pressure. Figure \ref{shear-setup} (a) sketches the problem domain and the displacement loading condition. The specimen has the same geometrical size as the one for Example 1 (Figure \ref{shear-setup} (a)) in \ref{example_1}. The bottom is fixed in all directions. Both lateral boundaries are prescribed with constant confining pressure. All fluid phase boundaries are no-flux. The specimen is discretized into 12,800 equal-sized material points with a characteristic length $\Delta x$ = 0.125 mm. The horizon $\delta$ is 0.25 mm. For the initial state, a uniform isotropic stress state of -100 kPa is prescribed on the skeleton, and apparent pre-consolidation pressure is -300 kPa. The initial suction is  20 kPa. Equivalently, the pore water pressure is -20 kPa. The constant confining pressure on the lateral boundaries is 84 kPa, which is imposed through an equivalent {total} force density \cite{song2018modeling}. For the loading protocol, the confining pressure is applied by a linear ramp function in 0.04 s. Then, the top boundary is subjected to a vertical compression with a rate of $3\times10^{-5}$ m/s. The total simulation time is 24 s.

\begin{figure}[h!]
\centering
    \includegraphics[width=0.35\textwidth]{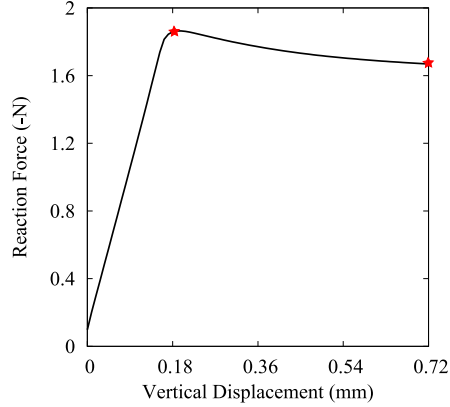}
  \caption{Reaction force versus the vertical displacement at the top boundary. }
    \label{plane-compr-forceplot}
\end{figure}%

\begin{figure}[h!]%
\centering
    \includegraphics[width=0.55\textwidth]{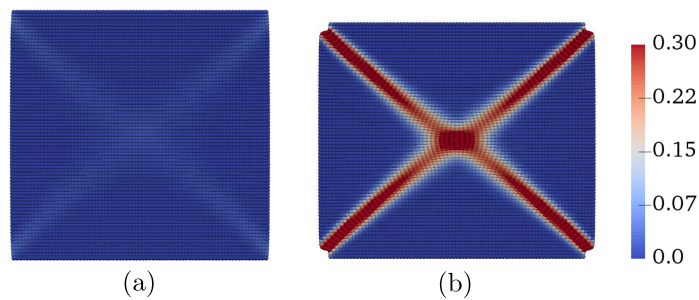}%
  \caption{Contours of equivalent plastic shear strain at displacements:(a) $u_y$ = 0.18 mm and (b) $u_y$ = 0.72 mm.}
    \label{ex2-shear-base}
\end{figure}

\begin{figure}[h!]%
\centering
    \includegraphics[width=0.55\textwidth]{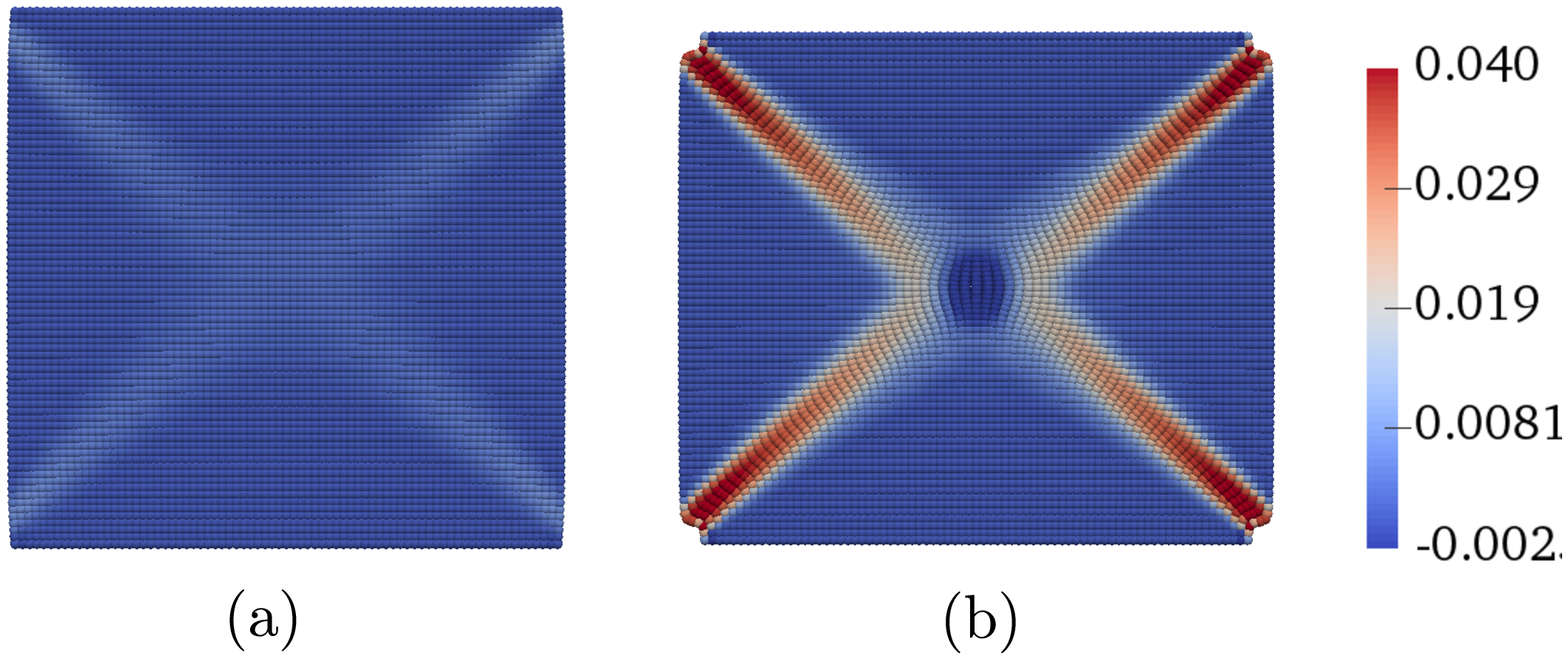}%
  \caption{Contours of plastic volume strain at displacements: (a) $u_y$ = 0.18 mm and (b) $u_y$ = 0.72 mm.}
  \label{ex2-volume-base}
\end{figure}

\begin{figure}[h!]%
\centering
    \includegraphics[width=0.55\textwidth]{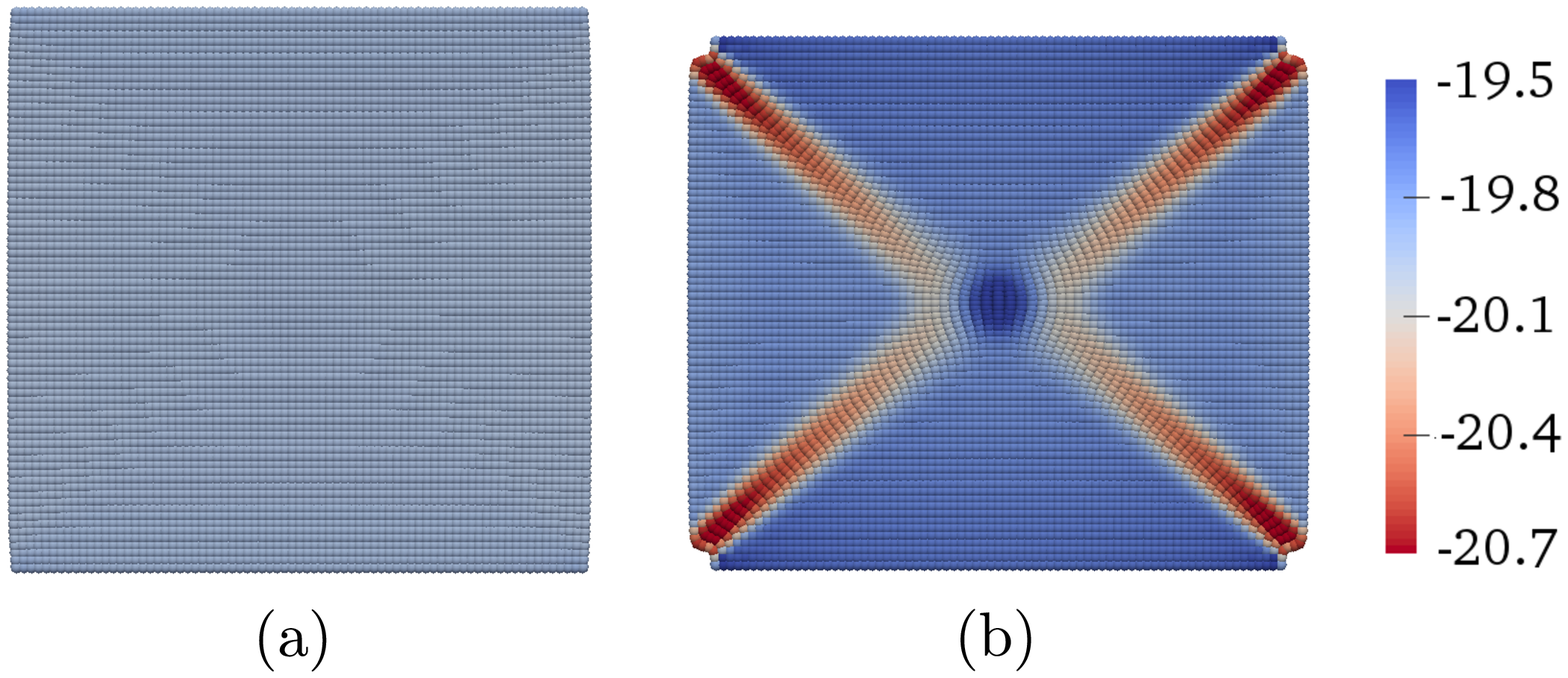}%
  \caption{Contours of pore water pressure (kPa) at displacements: (a) $u_y$ = 0.18 mm and (b) $u_y$ = 0.72 mm.}
    \label{ex2-pressure-base}
\end{figure}

\begin{figure}[h!]
\centering
    \includegraphics[width=0.55\textwidth]{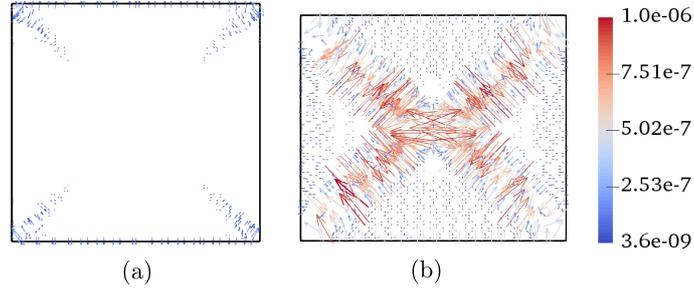}
  \caption{Contours of fluid flow vector and the magnitude at displacements: (a) $u_y$ = 0.18 mm and (b) $u_y$ = 0.72 mm.}
    \label{ex2-flux-base}
\end{figure}

Figure \ref{plane-compr-forceplot} plots the reaction force versus the vertical displacement of the top boundary. The markers in Figure \ref{plane-compr-forceplot} correspond to two applied displacements: (a) $u_y$  = 0.18 mm and (b) $u_y$ = 0.72 mm, respectively. The peak reaction force occurs at $u_y$ = 0.18 mm. Simulation results at both load stages are shown in {Figures \ref{ex2-shear-base}, \ref{ex2-volume-base}, \ref{ex2-pressure-base} and \ref{ex2-flux-base}}. {All contours depicted here are drawn on the deformed configuration of the specimen}. Figure \ref{ex2-shear-base} presents the snapshots of equivalent plastic shear strain. The results show the formation of two conjugate shear bands in the specimen. Figure \ref{ex2-volume-base} presents the snapshots of the plastic volumetric strains. The results in Figure \ref{ex2-volume-base} show that positive plastic volumetric strains occur in localized shear zones. Figure \ref{ex2-pressure-base} presents the snapshots of pore water pressure. Similar to the results for the deformation field in Figure \ref{ex2-shear-base} and \ref{ex2-volume-base}, two conjugate banded zones have formed in the pore water pressure field. Within the shear band, the pore water pressure is smaller than that outside the shear band. It is also evident that both the deformation band and the fluid pressure band have a similar finite thickness in that the same horizon is adopted for both the solid and the fluid. Figure \ref{ex2-flux-base} shows the snapshots of the fluid flux vector and the magnitude at the different load stages. The results show that the pore water flows into the shear band. This is due to the fact that the pore water within the shear band decreases (or equivalently suction increases) in the formation of dilative shear bands.

\subsubsection{Influence of horizon}

\begin{figure}[h!]
\centering
    \includegraphics[width=0.35\textwidth]{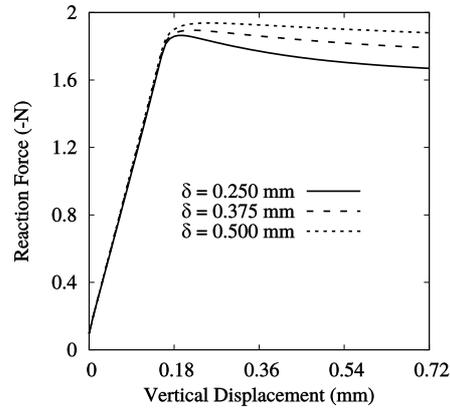}%
  \caption{Comparison of force versus applied vertical displacement for $\delta$  = 0.25 mm, (b) $\delta$ = 0.375 mm, and (c) $\delta$ = 0.5 mm.}
      \label{ex2-force-horizons}
\end{figure}

We study the influence of the horizon on the formation of shear bands under the symmetrical loading condition. In so doing, the previous simulation is repeated with two more horizons, $\delta$ = 0.375 mm and $\delta$ = 0.5 mm while keeping all other parameters unchanged. Figure \ref{ex2-force-horizons} compares the loading capacity curves (i.e., the applied loads versus the vertical displacement the top boundary) generated by the simulations with the three horizons, respectively. The results show that the loading response curves are slightly sensitive to the horizon size in the post-localization regime. The simulation with a larger horizon size generates a higher loading response at the same load step. 

\begin{figure}[h!]
\centering
    \includegraphics[width=0.6\textwidth]{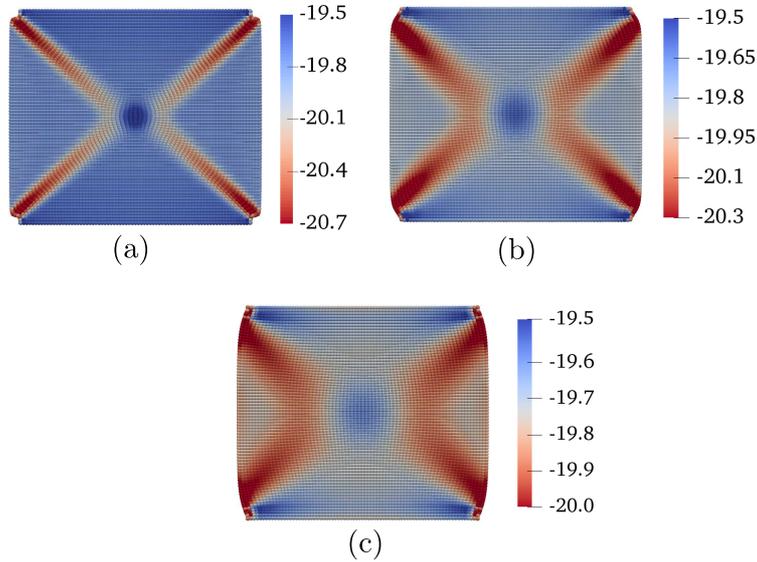}%
  \caption{Contours of pore water pressure at $u_y$ = 0.72 mm  for for (a) $\delta $ = 0.25 mm, (b) $\delta$ = 0.375 mm, and (c) $\delta$ = 0.5 mm.}
  \label{ex2-pressure-horizons}
\end{figure}

\begin{figure}[h!]
\centering
    \includegraphics[width=0.6\textwidth]{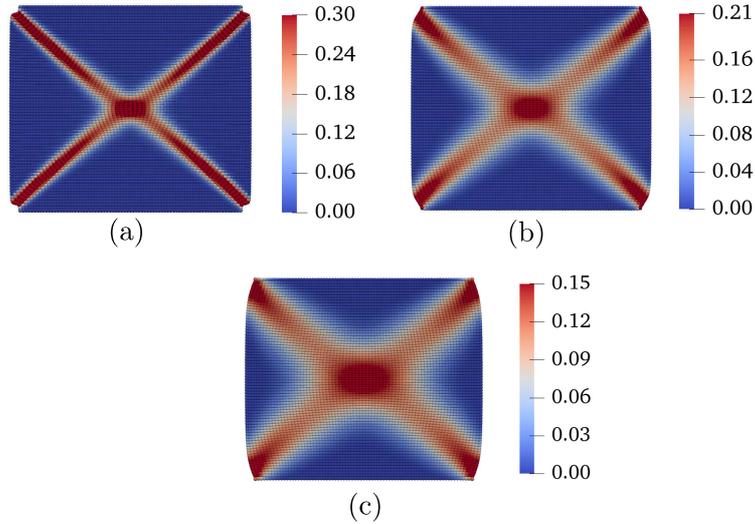}%
  \caption{Contours of equivalent plastic shear strain at $u_y$ = 0.72 mm for (a) $\delta$ = 0.25 mm, (b) $\delta$ = 0.375 mm, and (c) $\delta$ = 0.5 mm. }
  \label{ex2-shear-horizons}
\end{figure}

{Figures \ref{ex2-pressure-horizons} and \ref{ex2-shear-horizons}} show the contours of pore water pressure and plastic shear strain at the last load step of the simulations with the three horizons. Similar to the results of the shear example, the horizon size can influence the deformation and pore water pressure in the post-localization regime (e.g., the thickness of shear bands and the magnitude of the pore water pressure in the shear band). For instance, the simulation with the smallest horizon size generates the most intense plastic shear strain within the shear bands. 
 
 \subsubsection{Influence of hydraulic conductivity}

\begin{figure}[h!]%
\centering
    \includegraphics[width=0.35\textwidth]{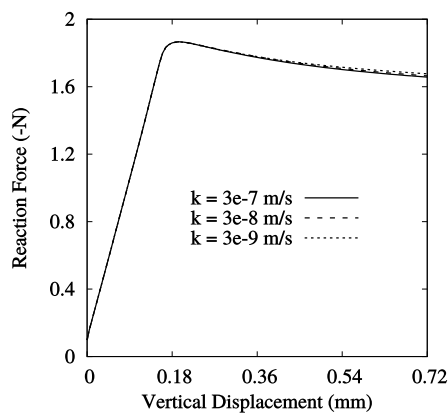}%
  \caption{Plot of reaction force versus vertical displacement for simulations with three intrinsic hydraulic conductivities: (a) $k$ = 3$\times 10^{-7}$ m/s, (b) $k$ = 3$\times 10^{-8}$ m/s, and (c) $k$ = 3$\times 10^{-9}$ m/s.}
  \label{ex2-force-permeability}
\end{figure}

\begin{figure}[h!]
\centering
    \includegraphics[width=0.6\textwidth,clip]{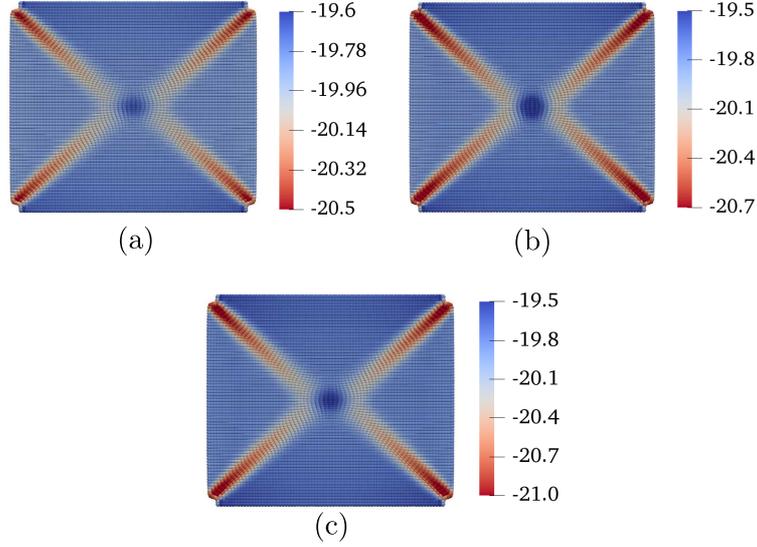}%
  \caption{Contours of pore water pressure (kPa) at $u_y$ = 0.72 mm for simulations with three intrinsic hydraulic conductivities: (a) $k$ = 3$\times 10^{-7}$ m/s, (b) $k$ = 3$\times 10^{-8}$ m/s, and (c) $k$ = 3$\times 10^{-9}$ m/s.}
  \label{ex2-pressure-permeability}
\end{figure}

\begin{figure}[h!]
\centering
    \includegraphics[width=0.6\textwidth,clip]{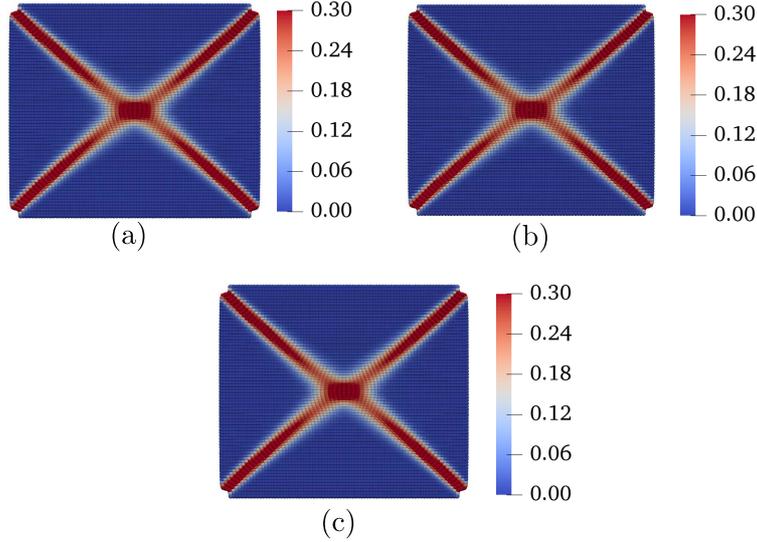}%
  \caption{Contours of equivalent plastic shear strain at $u_y$ = 0.72 mm for simulations with three intrinsic hydraulic conductivities: (a) $k$ = 3$\times 10^{-7}$ m/s, (b) $k$ = 3$\times 10^{-8}$ m/s, and (c) $k$ = 3$\times 10^{-9}$ m/s.}
  \label{ex2-shear-permeability}
\end{figure}

We analyze the influence of the intrinsic hydraulic conductivity on the numerical results. Assuming a constant viscosity (i.e., isothermal condition) and density (i.e., low water pressure) of pore water \cite{bear2013dynamics}, this analysis is equivalent to studying the influence of intrinsic permeability on the numerical results. For this purpose, we re-run the simulation in Section \ref{example_2}  with three different values of hydraulic conductivity: (a) $k$ = 3$\times 10^{-7}$ m/s, (b) $k$ = 3$\times 10^{-8}$ m/s, and (c) $k$ = 3$\times 10^{-9}$ m/s. Figure \ref{ex2-force-permeability} plots the reaction force versus the vertical displacement for the three simulations. Figure \ref{ex2-pressure-permeability} portrays the contour of pore water pressure at $u_y$= 0.72 mm for the three simulations. Figure \ref{ex2-shear-permeability} presents the contours of equivalent plastic shear strains at $u_y$= 0.72 mm.

The contours of water pressure show that the hydraulic conductivity affects the magnitude of the water pressure in the post-localization regime. For example, the results in Figure \ref{ex2-pressure-permeability} the simulation with a smaller value of intrinsic hydraulic conductivity generates a more considerable pressure change inside the shear band. Furthermore, the results in Figure \ref{ex2-shear-permeability} show that the intrinsic permeability has a negligible impact on the thickness of shear bands as well as the magnitude of equivalent plastic shear strain within the shear band. This observation may be due to the strong nonlocality for the fluid and solid phase in the coupled periporomechanics model. Note further study on the aspect shall be conducted in the future.

\subsubsection{Influence of spatial discretization}

\begin{figure}[h!]
\centering
    \includegraphics[width=0.35\textwidth]{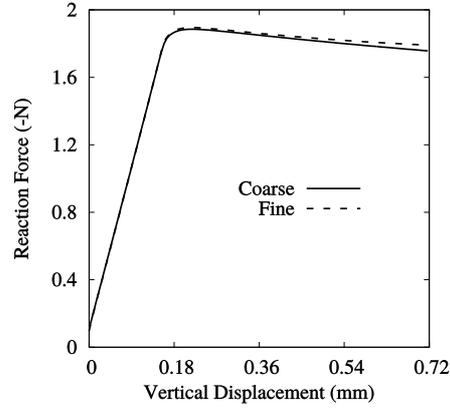}
  \caption{Plot of reaction force versus vertical displacement with two different spatial discretization schemes.}
  \label{ex2-force-meshes}
\end{figure}

To examine the influence of spatial discretization on the numerical results, we rerun the simulation in Section \ref{example_2} with two spatial discretization schemes: (a) the coarse one with 3200 material points ($\Delta x$ = 0.25 mm) and (b) the fine one with 12,800 material points ($\Delta x$ = 0.125 mm), respectively. To compare the results, the same horizon $\delta$ = 0.375 mm is adopted for both simulations. Figure \ref{ex2-force-meshes} plots the reaction force over the vertical displacement at the top boundary for the simulations with two spatial discretization schemes. It is apparent that the force predicted in both simulations is almost identical.

\begin{figure}[h!]
\centering
    \includegraphics[width=0.5\textwidth]{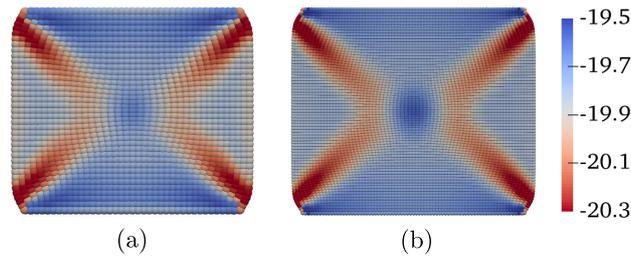}%
  \caption{Contours of pore water pressure (kPa) at $u_y$ = 0.72 mm using (a) coarse discretization and (b) fine discretization in spatial domain.}
  \label{ex2-pressure-meshes}
\end{figure}%

\begin{figure}[h!]
\centering
    \includegraphics[width=0.5\textwidth]{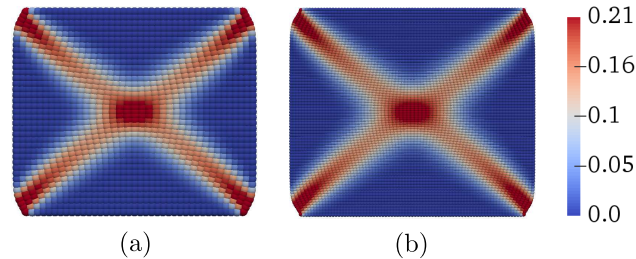}%
  \caption{Contours of equivalent plastic shear strains at $u_y$ = 0.72 mm using (a) coarse discretization and (b) fine discretization in spatial domain.}
  \label{ex2-shear-meshes}
\end{figure}

Figure \ref{ex2-pressure-meshes} depicts the contours of water pressure superimposed on the deformed configuration at $u_y$ = 0.72 mm. Figure \ref{ex2-shear-meshes} plots the contours of equivalent plastic shear strains at $u_y$= 0.72 mm. The results in both figures evidence that the equivalent plastic shear strain and water pressure are insensitive to spatial discretizations. For both cases, the plastic shear deformation and the water pressure localize in a well-defined zone of finite thickness that is correlated to the horizon. It may be concluded that the proposed coupled periporomechanics is not sensitive to the spatial discretization, which is different from the standard finite element method. Thus, the proposed periporomechanics can be a valuable tool in the analysis of localization phenomena in unsaturated porous media.

 \subsection{Example 3: Bearing capacity of a footing on unsaturated soils} \label{example_3}

\begin{figure}[h!]
\centering
    \includegraphics[width=0.5\textwidth]{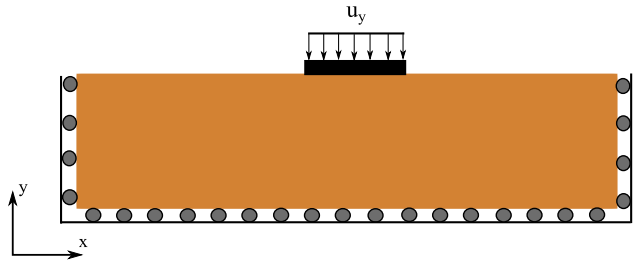}
  \caption{Problem setup.}
  \label{example3-strip-setup}
\end{figure}

This example deals with the bearing capacity of footing on unsaturated soils. Specifically, this example aims to show the proposed computational periporomechanics model can qualitatively simulate the localized failure mode of unsaturated soils under a footing. Figure \ref{example3-strip-setup} sketches the geometry of the problem domain and the applied boundary conditions. A footing on which the displacement load is imposed is placed on the top of the soil domain. The dimension of the soil domain is 4 m $\times$ 1 m $\times$ 0.02 m. The footing's dimension is 0.8 $\times$ 0.01 $\times$ 0.02 m. As shown in Figure \ref{example3-strip-setup}, the lateral boundaries are constrained from horizontal deformation and are free to deform vertically. The bottom is constrained from vertical deformation. All fluid phase boundaries are no-flux conditions. The problem domain is discretized into 60,000 mixed material points. The multiphase length scale of $\delta$ is assumed 0.022 m. 

The material parameters for the skeleton are as follows. Bulk modulus is 55.56 $\times 10^3$ kPa, the shear modulus is 18.52 $\times 10^3$ kPa, swelling index $\kappa$ is 0.05, the slope of critical state line slope is 1.0, compression index is 0.15, and specific volume at unit pre-consolidation pressure is 2.0, initial apparent pre-consolidation pressure is -200 kPa, and other parameters are the same adopted in examples 1 and 2. The parameters for soil-water characteristic curve are: $S_1$ = 0.0, $S_2$ = 1.0, $s_a$ = 10 kPa, $n$ = 1.25, and $m$ = 0.2. The soil density is 2500 kg/$\mathrm{m}^3$, water density is 1000 kg/$\mathrm{m}^3$, initial porosity is 0.25, and hydraulic conductivity is 1$\times 10^{-7}$ m/s. The initial suction in the soil is assumed 20 kPa. The corresponding degree of saturation and relative permeability are 0.8 and 0.01, respectively. The loading protocol is divided into two steps. In step 1, the soil's initial geostatic stress state is prescribed by imposing a gravity load on the soil. Then, a null deformation condition is assigned to the soil before a load is imposed on the footing. In step 2, the footing is given a vertical displacement $u_y$ = 0.04 m at the constant rate of 1.67 mm/s in 15 s.

\begin{figure}[h!]
\centering
    \includegraphics[width=0.5\textwidth]{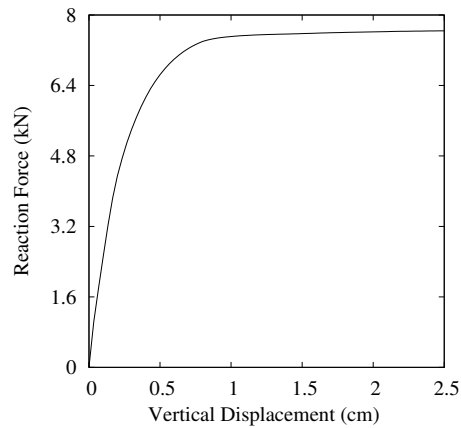}
  \caption{Plot of the reaction force on the footing versus the applied vertical displacement.}
  \label{example3-reaction}
\end{figure}

 \begin{figure}[h!]
\centering
    \includegraphics[width=0.5\textwidth]{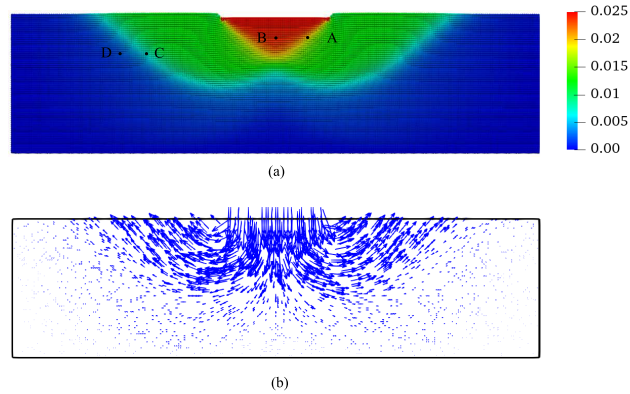}
  \caption{Contour of (a) the magnitude and (b) direction of the displacement vector of the soil at $u_y$ = 2.5 cm.}
  \label{example3-displacement-contours}
\end{figure} 

 \begin{figure}[h!]
\centering
    \includegraphics[width=0.5\textwidth]{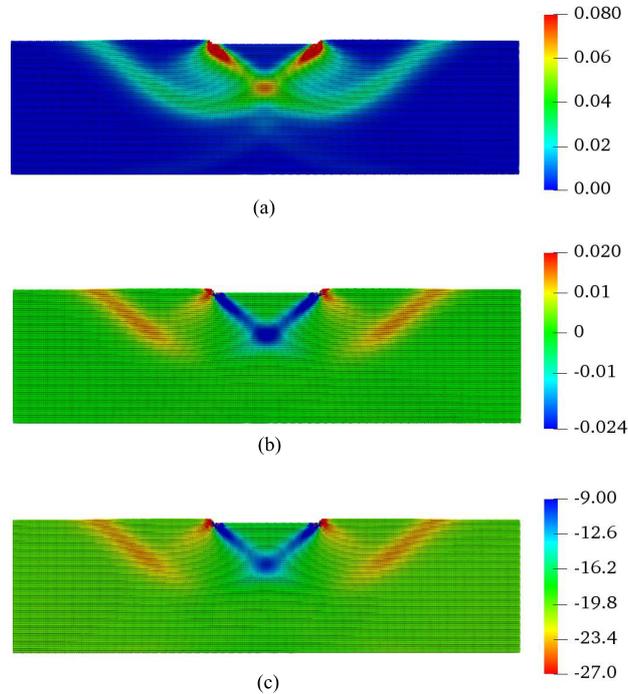}
  \caption{Contours of (a) equivalent plastic shear strain, (b) plastic volumetric strain, and (c) pore water pressure (kPa) at $u_y $= 2.5 cm.}
  \label{example3-base-results}
\end{figure} 

\begin{figure}[h!]
\centering
    \includegraphics[width=0.8\textwidth]{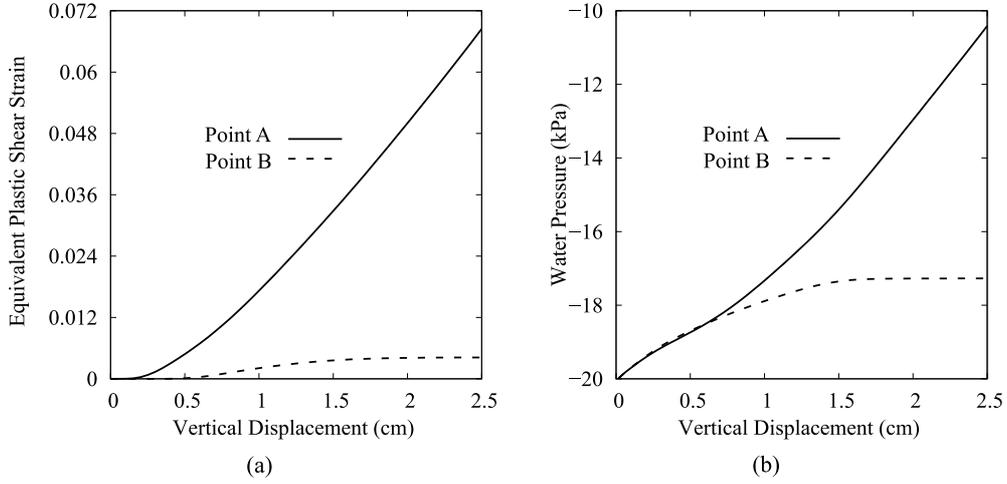}
  \caption{Variation of (a) equivalent plastic shear strain and (b) pore water pressure versus the footing's vertical displacement at the points A and B (inside and outside banded compaction zones, respectively).}
  \label{example3-inside-outside-compact}
\end{figure}

\begin{figure}[h!]
\centering
    \includegraphics[width=0.8\textwidth]{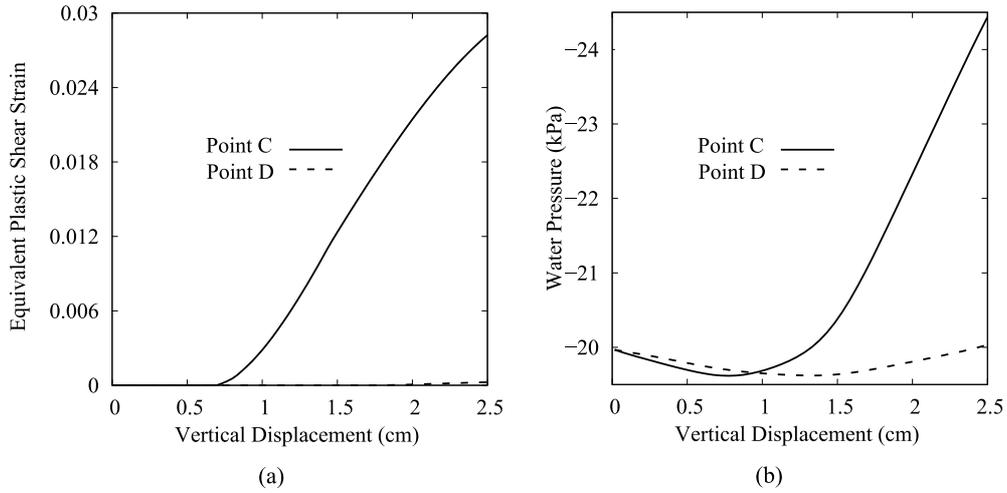}
  \caption{Variation of (a) equivalent plastic shear strain and (b) pore water pressure at the points C and D versus the footing's vertical displacement (inside and outside banded dilation zones, respectively).}
  \label{example3-inside-outside-dilatant}
\end{figure}

Figure \ref{example3-reaction} plots the reaction force of the footing versus the vertical displacement. As shown in Figure \ref{example3-reaction}, the bearing capacity of the footing is around 7.5 kN. The contour of the magnitude of displacement is presented in Figure \ref{example3-displacement-contours} (a). The displacement vector is plotted in Figure \ref{example3-displacement-contours} (b) in which the arrowhead denotes the direction, and the length of the arrow tail shows the scaled magnitude. The contours of equivalent plastic shear strain, plastic volumetric strain (or plastic dilatation), pore water pressure in the soil at $u_y$ = 2.5 cm are depicted in Figure \ref{example3-base-results} (a), (b) and (c). The results show that a typical localized shear failure mode of soil under shallow foundations \cite{terzaghi1951theoretical}. For instance, the settlement of the footing has caused the formation of symmetric slip surfaces in the soil and soil upheaval on both sides of the footing. It is apparent from Figure \ref{example3-base-results} (b) that the shear bands under the footing are compactive while those associated with soil upheaval are dilatant. It is also revealed from Figure \ref{example3-base-results} (c) that matric suction decreases (pore pressure rises) in the compactive shear bands beneath the footing while matric suction increases in the dilative shear bands. The results are corroborated by the data shown in Figures \ref{example3-inside-outside-compact} and \ref{example3-inside-outside-dilatant} which present the variation of the equivalent plastic strain and pore water pressure versus the footing's vertical displacement at points inside and outside the shear bands denoted in Figure \ref{example3-displacement-contours} (a).

 \begin{figure}[h!]
\centering
    \includegraphics[width=0.4\textwidth]{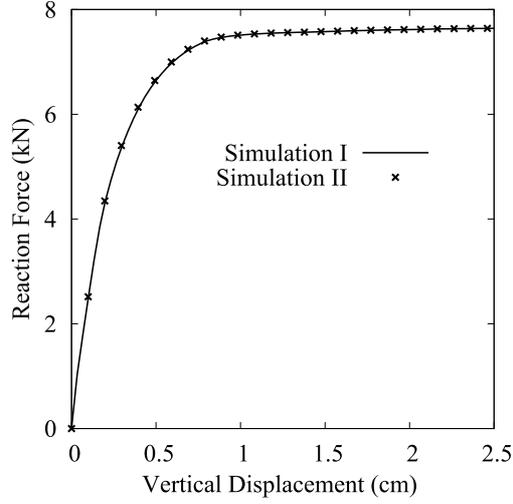}
  \caption{Plot of reaction force on the footing over the applied vertical displacement for Simulation I and Simulation II.}
  \label{example3-reaction-discretization}
\end{figure}

 \begin{figure}[h!]
\centering
    \includegraphics[width=0.5\textwidth]{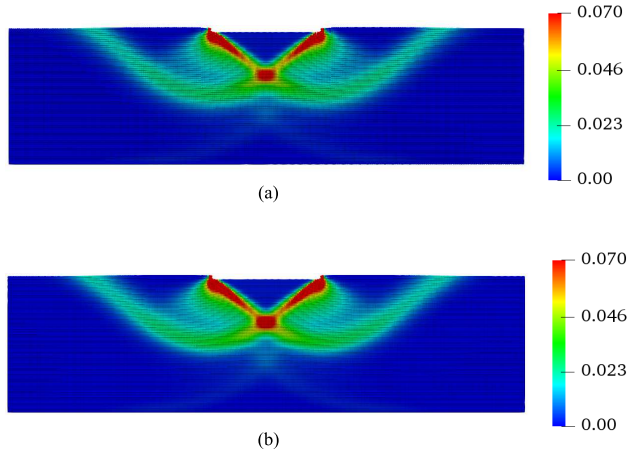}
  \caption{Contours of equivalent plastic shear strain for (a) Simulation I and (b) Simulation II at $u_y$ = 2.5 cm.}
  \label{example3-shear-comparison}
\end{figure}

 \begin{figure}[h!]
\centering
    \includegraphics[width=0.5\textwidth]{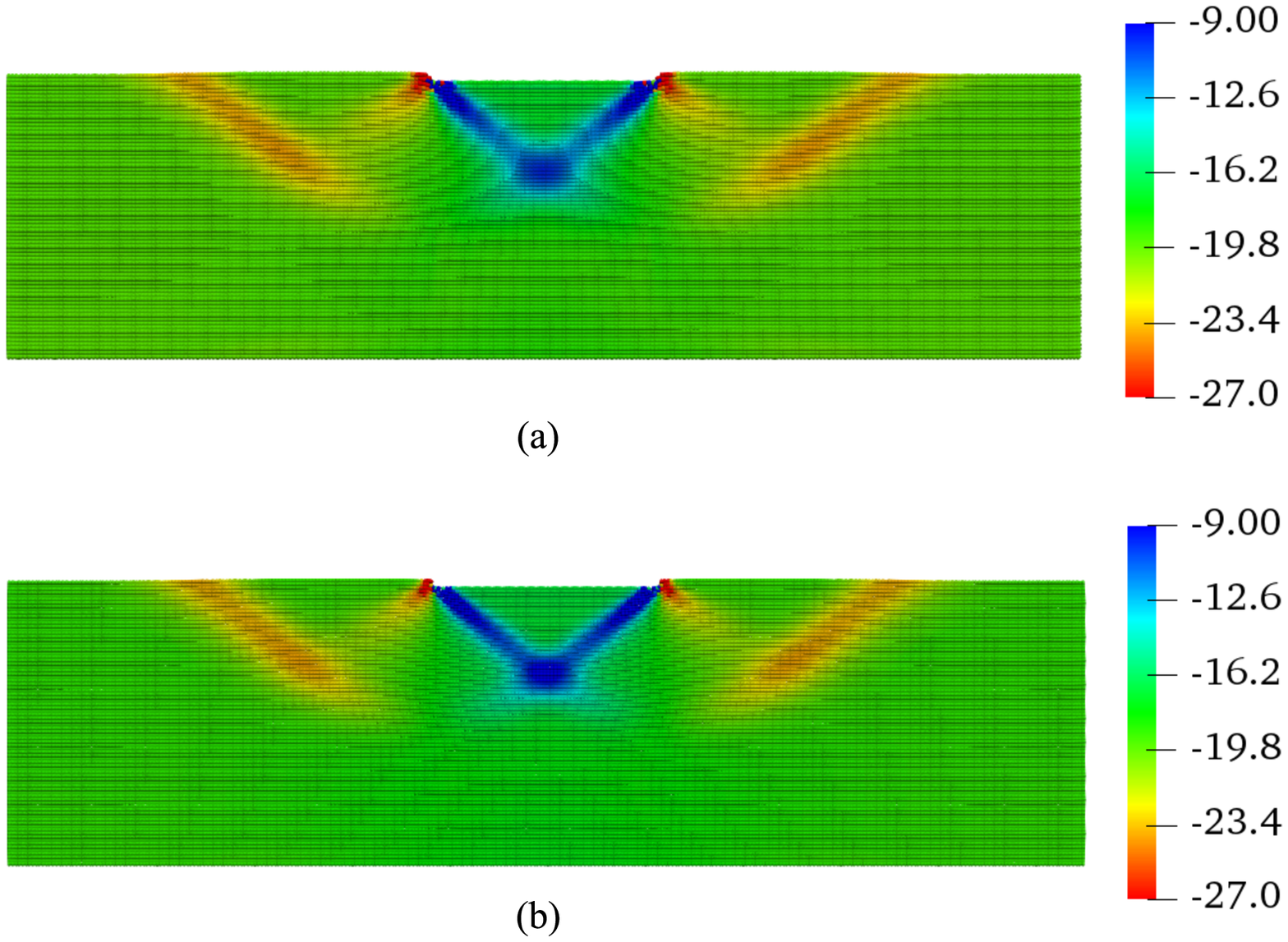}
  \caption{Contours of pore water pressure for (a) Simulation I and (b) Simulation II at $u_y$ = 2.5 cm.}
  \label{example3-pressure-comparison}
\end{figure}

To test the solution's uniqueness, we rerun the simulation with a relatively finer spatial discretization consisting of 100,000 mixed material points. The same length scale and input parameters are adopted. The comparison between the two simulations is presented in Figures \ref{example3-reaction-discretization}, \ref{example3-shear-comparison} and \ref{example3-pressure-comparison}. Here simulations I and II stand for the simulations with 60,000 and 100,000 material points respectively. Figure \ref{example3-reaction-discretization} shows that the reaction force curves from both simulations are almost identical. Moreover, as shown in Figures \ref{example3-shear-comparison} and \ref{example3-pressure-comparison} respectively the contours of equivalent plastic shear strain and pore water pressure at $u_y$ = 2.5 cm from both simulations are in good agreement. It is noted that all the numerical simulations are run in parallel on a supercomputer. To show the performance of the Open MPI implementation and scalability of the numerical code, the total wall-clock time for both simulations run in parallel with different number of CPU (central processing unit) processors is summarized in Table \ref{wall-clock-time-comp}. With the same number of CPU processors, the wall-clock time of simulation I with 60,000 multiphase material points is much shorter than that of simulation II with 100,000 multiphase material points. Doubling the number of CPU processors reduces the wall-clock time dramatically for both simulations, although the reduction of time is less than half. The latter can be because more message passing occurs between processors for the simulations with 64 CPUs than those with 32 CPUs.

\begin{table}[h!]
\centering
\caption{Comparison of wall-clock time for simulations I and II}
\label{wall-clock-time-comp}    
\begin{tabular}{c|cc}
\noalign{\smallskip}
\noalign{\smallskip}\hline\noalign{\smallskip}
Number of CPU Processors &\multicolumn{2}{c}{Wall-clock Time} \\
\cline{2-3}
\\[-10pt]
 & Simulation I (60,000 material points) & Simulation II (100,000 material points)\\
\hline
\\[-5pt]
32 & 15,600 s &  44,500 s \\
64 & 8,900  s &  33,000 s\\
\noalign{\smallskip}\hline
\end{tabular}
\end{table}

\section{Closure}

In this article, we have implemented a computational periporomechanics model for simulating localized failure in unsaturated porous media. The mathematical model is based on the peridynamic state concept and the effective force state concept. It is assumed that the pore air pressure in unsaturated porous media is passive (i.e., atmospheric pressure). The coupled governing equations are integral-differential equations without assuming the continuity of solid displacement and pore water pressure (or matric suction). Both equations have an intrinsic two-phase length scale - the horizon. In the coupled periporomechanics model, it is assumed that unsaturated porous media consist of equal-sized mixed material points, which have two kinds of degrees of freedom, i.e., the displacement and pore water pressure. In line with this hypothesis, we have proposed an implicit multiphase Lagrangian meshless method to numerically implement the periporomechanics model. The numerical implementation is validated by simulating Terzaghi's consolidation problem, drainage of a soil column, and Mandel's slab problem. Numerical examples under different scenarios are conducted to demonstrate that the computational periporomechanics model is robust in modeling the formation of strain localization in unsaturated porous media. It is shown that the intrinsic multiphase length scale can be utilized to characterize the finite thickness of the localization of deformation and pore water pressure. Length scale sensitivity analysis has shown that the multiphase length scale can impact the formation of shear bands and the loading capacity of unsaturated soils in the post-localization regime. It was found from the numerical analysis that spatial discretization and intrinsic hydraulic conductivity have mild effects on localized failure of unsaturated soils. This may be due to the fact that the identical length scale (i.e., the horizon) is chosen for both the pore fluid and solid phases. It is shown from the footing example that the computational periporomechanics model can replicate the shear failure modes of soils under a strip footing at the field scale. 

\section*{Acknowledgment}
Support for this work was provided by the US National Science Foundation under Contract numbers 1659932 and 1944009. The authors are grateful to the anonymous reviewers for their constructive reviews.

\bibliography{ijnme_cpm_2020}

\end{document}